# Bounds for the covariance of functions of infinite variance stable random variables with applications to central limit theorems and wavelet-based estimation

VLADAS PIPIRAS[1], MURAD S. TAQQU[2] and PATRICE ABRY[3]

[1]*Department of Statistics and Operations Research, UNC-CH, CB #3260, Chapel Hill, NC 27599-3260, USA. E-mail: pipiras@email.unc.edu*

[2]*Department of Mathematics and Statistics, Boston University, 111 Cummington St., Boston, MA 02215, USA. E-mail: murad@math.bu.edu*

[3]*CNRS UMR 5672, Ecole Normale Supérieure de Lyon, Laboratoire de Physique, 69 364 Lyon Cedex 07, France. E-mail: pabry@physique.ens-lyon.fr*

We establish bounds for the covariance of a large class of functions of infinite variance stable random variables, including unbounded functions such as the power function and the logarithm. These bounds involve measures of dependence between the stable variables, some of which are new. The bounds are also used to deduce the central limit theorem for unbounded functions of stable moving average time series. This result extends the earlier results of Tailen Hsing and the authors on central limit theorems for bounded functions of stable moving averages. It can be used to show asymptotic normality of wavelet-based estimators of the self-similarity parameter in fractional stable motions.

*Keywords:* central limit theorem; covariance; dependence measures; linear fractional stable motion; moving averages; self-similarity parameter estimators; stable distributions; wavelets

## 1. Introduction

Consider a (non-Gaussian) symmetric $\alpha$-stable ($S\alpha S$, in short), $\alpha \in (0, 2)$, two-dimensional random vector $(\xi, \eta)$. By definition (see, e.g., Samorodnitsky and Taqqu [19]), this means that any linear combination $u\xi + v\eta$ is an $S\alpha S$ random variable, that is, it has a characteristic function $Ee^{i\theta(u\xi+v\eta)} = e^{-\sigma^\alpha |\theta|^\alpha}$, where $\sigma > 0$ is the scale coefficient, denoted by $\|u\xi + v\eta\|_\alpha$. We will sometimes assume that $\alpha > 1$ so that the variables $\xi$ and $\eta$ have finite mean. These random variables, however, always have infinite variance since $\alpha < 2$. Nevertheless, functions $K(\xi)$ and $L(\eta)$ of these random variables may have finite variance, as is the case for bounded functions $K, L$ or unbounded functions $K(x) = |x|^{\beta/2}$ with $\beta < \alpha/2$, $L(y) = \ln|y|$. For such functions, the covariance

$$\mathrm{Cov}(K(\xi), L(\eta)) \tag{1.1}$$

is well defined.

We are interested here in bounds on the covariance (1.1) when $\xi$ and $\eta$ are dependent. (Examples of such $\xi, \eta$ are provided in Section 2.) We seek bounds on (1.1) in terms of workable measures of dependence between $\xi$ and $\eta$. Since $\xi$ and $\eta$ have infinite variance, dependence measures other than the usual covariance must be used. Some alternative measures of dependence are available in the stable case, such as codifference and covariation (Samorodnitsky and Taqqu [19]). We will use a measure similar to the covariation and also introduce a new measure of dependence.

Bounds on the covariance (1.1) are useful for several reasons. From a general perspective, they complement the work on bounds of covariances of other random variables. In some cases, for example, Gaussian variables and integer power functions, these covariances can be computed explicitly in terms of covariance between the two variables. In other cases, more general dependence measures, such as various types of mixing, are studied. As with bounds on covariances for other variables, bounds on (1.1) are useful for establishing limit results. In fact, we will use the obtained bounds to deduce the central limit theorem (CLT) for partial sums

$$\frac{1}{N^{1/2}} \sum_{n=1}^{N} (K(\xi_n) - EK(\xi_n)), \tag{1.2}$$

where $K$ is as in (1.1) and $\xi_n$ is an $S\alpha S$ moving average time series with $\alpha \in (1, 2)$. As $K$ can be unbounded, this result extends the CLTs for partial sums (1.2), with bounded functions $K$ considered in Hsing [11] and Pipiras and Taqqu [16]. A related paper of Wu [23], which allows functions $K$ to be unbounded, will be compared with our result below. We also apply the obtained CLT to establish the asymptotic normality of wavelet-based estimators of the self-similarity parameter in linear fractional stable motion. The obtained asymptotic normality result is, to our knowledge, the first of its type.

The paper is organized as follows. Section 2 contains some preliminaries and statements of the main results on bounds of the covariance (1.1). A useful measure of dependence is studied in Section 3. The proofs of the main results on covariance bounds can be found in Section 4. Section 5 contains auxiliary results. In Section 6, we show the CLT for partial sums (1.2) and in Section 7, we apply it to deduce asymptotic normality of wavelet-based estimators in linear fractional stable motion.

## 2. Results on covariance bounds

Let $(\xi, \eta)$ be an $S\alpha S$, $\alpha \in (0, 2)$, random vector, as in Section 1. We will be concerned throughout only with distributional properties of $(\xi, \eta)$. Therefore (Samorodnitsky and

Taqqu [19]), we may suppose, without loss of generality, that

$$(\xi, \eta) = \left( \int_S f(s) M(\mathrm{d}s), \int_S g(s) M(\mathrm{d}s) \right), \tag{2.1}$$

where $(S, \mu)$ is a measure space, $f, g \in L^\alpha(S, \mu)$ and $M$ is the so-called $S\alpha S$ random measure on $S$ with the control measure $\mu$. Heuristically, $M(\mathrm{d}s)$, $s \in S$, can be viewed as a sequence of independent $S\alpha S$ random variables with scale coefficients $\mu(\mathrm{d}s)$. More precisely, for any Borel set $A$, $M(A)$ is an $S\alpha S$ random variable with scale coefficient $\|M(A)\|_\alpha = \mu(A)^{1/\alpha}$. Moreover, if $A_1, \ldots, A_n$ are disjoint Borel sets, $M(A_1 \cup \cdots \cup A_n) = M(A_1) + \cdots + M(A_n)$ and $M(A_1), \ldots, M(A_n)$ are independent. Because the random variables $M(A_1), \ldots, M(A_n)$ are $S\alpha S$, the resulting vector $(\xi, \eta)$ is $S\alpha S$. (For more details, see Chapter 3 in Samorodnitsky and Taqqu [19].) With the notation $\|f\|_\alpha = (\int_S |f(s)|^\alpha \mu(\mathrm{d}s))^{1/\alpha}$, the relation (2.1) also means that the characteristic function of $(\xi, \eta)$ is given by

$$E \mathrm{e}^{\mathrm{i}(u\xi + v\eta)} = \mathrm{e}^{-\|uf + vg\|_\alpha^\alpha}, \tag{2.2}$$

that is, $\|u\xi + v\eta\|_\alpha = \|uf + vg\|_\alpha$.

According to Samorodnitsky and Taqqu [19], the space $(S, \mu)$ in (2.1) can be taken as $S = (0, 1)$, $\mu =$ Lebesgue measure. Other choices are possible and are often used.

**Example 2.1.** Consider $S\alpha S$ random variables

$$\xi = \sum_{k=-\infty}^{\infty} f_k \varepsilon_k, \qquad \eta = \sum_{k=-\infty}^{\infty} g_k \varepsilon_k,$$

where $\sum |f_k|^\alpha < \infty$, $\sum |g_k|^\alpha < \infty$ and $\{\varepsilon_n\}$ is an i.i.d. sequence of $S\alpha S$ random variables (with scaling coefficient 1). Then, $S = \mathbb{Z}$, $\mu =$ counting measure, $f(n) = f_n$ and $g(n) = g_n$.

**Example 2.2.** Consider an $S\alpha S$ stationary time series given by the moving average representation

$$\xi_n = \int_\mathbb{R} a(n - x) M(\mathrm{d}x), \tag{2.3}$$

where $a \in L^\alpha(\mathbb{R}, \mathrm{d}x)$ and $M$ is an $S\alpha S$ random measure on $\mathbb{R}$ having the Lebesgue control $\mathrm{d}x$. Take $\xi = \xi_0$ and $\eta = \xi_n$. Then, $(S, \mu) = (\mathbb{R}, \mathrm{d}x)$ and $f(\cdot) = a(-\cdot), g(\cdot) = a(n - \cdot)$.

We shall use the following quantities related to $\xi$ and $\eta$.

**Definition 2.1.** Set

$$[\xi, \eta]_1^* = \int_S |f(s)|^{\alpha - 1} |g(s)| \mu(\mathrm{d}s), \qquad [\xi, \eta]_1 = [\xi, \eta]_1^* + [\eta, \xi]_1^*, \tag{2.4}$$

$$[\xi, \eta]_2 = \int_S |f(s) g(s)|^{\alpha/2} \mu(\mathrm{d}s). \tag{2.5}$$

Observe that, unlike $[\xi,\eta]_1^*$, the quantities $[\xi,\eta]_1$ and $[\xi,\eta]_2$ are symmetric in their arguments $\xi$ and $\eta$. Observe also that, by Hölder's inequality, $[\xi,\eta]_2 < \infty$ with $\alpha \in (0,2)$ and $[\xi,\eta]_1^* < \infty$ with $\alpha \in (1,2)$. The quantities $[\xi,\eta]_1$ and $[\xi,\eta]_2$ can be viewed as measures of dependence between the random variables $\xi$ and $\eta$. In particular, the dependence measure $[\xi,\eta]_1$ is related to the covariation measure of a stable random vector $(\xi,\eta)$. The covariation, however, does not involve absolute values (see Samorodnitsky and Taqqu [19], Section 2.7). Observe, also, that the representation (2.1) of random variables $\xi$ and $\eta$ is not unique: for example, another representation could be obtained from (2.1) by a change of variables. It is therefore important to note that the quantities in (2.4)–(2.5) are invariant with respect to the representation (2.1) of $\xi$ and $\eta$, as stated in the following proposition. See Section 4 for its proof.

**Proposition 2.1.** *The quantities $[\xi,\eta]_1^*$, $[\xi,\eta]_1$ and $[\xi,\eta]_2$ in (2.4)–(2.5) do not depend on the representation (2.1) of $(\xi,\eta)$.*

We shall often assume that $\xi$ and $\eta$ satisfy the following conditions.

**Assumption (A1).** *There is a constant $\varepsilon_1 > 0$ such that*

$$\|f\|_\alpha^{\alpha/2}\|g\|_\alpha^{\alpha/2} - \int_S |f(s)g(s)|^{\alpha/2}\mu(\mathrm{d}s) \geq \varepsilon_1 \|f\|_\alpha^{\alpha/2}\|g\|_\alpha^{\alpha/2}. \tag{2.6}$$

**Assumption (A2).** *There is a constant $\varepsilon_2 > 0$ such that, for all $u,v \in \mathbb{R}$,*

$$\|uf + vg\|_\alpha^\alpha \geq \varepsilon_2(\|uf\|_\alpha^\alpha + \|vg\|_\alpha^\alpha). \tag{2.7}$$

For example, the variables $\xi$ and $\eta$ in Example 2.2 satisfy these assumptions for sufficiently large $n$, where the constants $\varepsilon_1$ and $\varepsilon_2$ do not depend on $n$. Indeed, as proved in Section 4, we have the following.

**Proposition 2.2.** *Let $\alpha \in (0,2)$ and $\{\xi_n\}_{n\in\mathbb{Z}}$ be an $S\alpha S$ average time series from Example 2.2. Then, for sufficiently large $n$, the random variables $\xi = \xi_0$ and $\eta = \xi_n$ satisfy Assumptions* (A1) *and* (A2), *where the constants $\varepsilon_1$ and $\varepsilon_2$ do not depend on $n$.*

We shall now state the assumptions on the functions $K, L$ in (1.1) which are used in this work. We shall assume the following.

**Assumption (K1).** *There are $\beta \in (0,\alpha/2)$ and $x_0 > 0$ such that, for $|x| > x_0$,*

$$|K(x)| \leq \mathrm{const}|x|^\beta. \tag{2.8}$$

**Assumption (K2).** *The function $K(x)/x$ is non-increasing for $x > x_0$ and non-decreasing for $x < -x_0$.*

**Assumption (K3).** *The function $K$ is bounded on $|x| \leq x_0$ except possibly at a finite number of points $x_1, \ldots, x_p$, around which the function $K$ is integrable.*

Assumption (K1) is quite natural since it ensures that $EK(\xi_1)^2 < \infty$. Assumption (K2) is rather technical and depends on our method of proof. Examples of functions $K$ satisfying (K1)–(K3) are $|x|^\beta$ with $\beta \in (-1, \alpha/2)$, $(\ln |x|)^m$ with $m \geq 1$ and many others. The choice of assumptions is also motivated by the applications to central limit theorems in Sections 6 and 7.

The following results, proved in Section 4, provide bounds for the covariance $\text{Cov}(K(\xi), L(\eta))$ under various assumptions on the functions $K$ and $L$. In Theorem 2.1, we suppose that $K$ and $L$ are any two integrable functions. Theorem 2.2 concerns the functions $K$ and $L$ that are non-zero only for large values of the argument and which satisfy Assumption (K2). Theorem 2.3 is a combination of the previous two results. Finally, Theorem 2.4 is a consequence of Theorems 2.1, 2.2 and 2.3.

**Theorem 2.1.** *Let $\alpha \in (0, 2)$ and $\xi, \eta$ be two $S\alpha S$ random variables (2.1) that satisfy Assumption (A1). Let also $K, L \in L^1(\mathbb{R})$. Then,*

$$|\text{Cov}(K(\xi), L(\eta))| \leq C \|K\|_1 \|L\|_1 [\xi, \eta]_2, \qquad (2.9)$$

*where the constant $C$ depends on $\varepsilon_1$, $\alpha$, $\|\xi\|_\alpha$ and $\|\eta\|_\alpha$.*

**Theorem 2.2.** *Let $\alpha \in (1, 2)$ and $\xi, \eta$ be two $S\alpha S$ random variables (2.1) that satisfy Assumptions (A1) and (A2). Also, let $K$ and $L$ be two functions satisfying Assumptions (K1) and (K2) with $x_0 = b \geq 1$. Set $K_b(x) = K(x) 1_{\{|x| > b\}}$ and $L_b(x) = L(x) 1_{\{|x| > b\}}$, $x \in \mathbb{R}$. Then*

$$|\text{Cov}(K_b(\xi), L_b(\eta))| \leq C b^{2\beta - \alpha} ([\xi, \eta]_1 + [\xi, \eta]_2), \qquad (2.10)$$

*where the constant $C$ depends on $\varepsilon_1$, $\varepsilon_2$, $\alpha$, $\beta$, $\|\xi\|_\alpha$ and $\|\eta\|_\alpha$.*

**Theorem 2.3.** *Let $\alpha \in (1, 2)$ and $\xi, \eta$ be two $S\alpha S$ random variables (2.1) that satisfy Assumptions (A1) and (A2). Also, let $K_b$ be a function as in Theorem 2.2 and $L \in L^1(\mathbb{R})$. Then,*

$$|\text{Cov}(K_b(\xi), L(\eta))| \leq C b^{\beta - 1} \log(b + 1) \|L\|_1 ([\eta, \xi]_1^* + [\eta, \xi]_2), \qquad (2.11)$$

*where the constant $C$ depends on $\varepsilon_1$, $\varepsilon_2$, $\alpha$, $\beta$, $\|\xi\|_\alpha$ and $\|\eta\|_\alpha$.*

**Theorem 2.4.** *Let $\alpha \in (1, 2)$ and $\xi, \eta$ be two $S\alpha S$ random variables (2.1) that satisfy Assumptions (A1) and (A2). Also, let $K$ and $L$ be two functions that satisfy Assumptions (K1)–(K3). Then*

$$|\text{Cov}(K(\xi), L(\eta))| \leq C([\xi, \eta]_1 + [\xi, \eta]_2), \qquad (2.12)$$

*where the constant $C$ depends on $\varepsilon_1$, $\varepsilon_2$, $K$, $L$, $\|\xi\|_\alpha$ and $\|\eta\|_\alpha$.*

**Remark 2.1.** Observe that the bounds (2.10) and (2.11) involve the cut-off parameter $b$. In particular, as $b \to \infty$, both bounds converge to 0, which is consistent with the fact that the two covariances also tend to 0. These results are therefore acceptable in the sense that the effects of the cut-off $b$ and the dependence are separated in the bounds. We also note that explicit expressions for the constants $C$ in the bounds of the covariances above can be deduced from the proofs of the results, but these are not pursued here.

## 3. Measure of dependence

In this section, we establish some results on a measure of dependence $U_{\xi,\eta}$ between $\xi$ and $\eta$ defined below. These results will be used in Section 4.

**Definition 3.1.** *For $S\alpha S$ random variables $\xi$ and $\eta$ in (2.1) and $u, v \in \mathbb{R}$, set*

$$U_{\xi,\eta}(u,v) = E e^{iu\xi + iv\eta} - E e^{iu\xi} E e^{iv\eta} = e^{-\|uf+vg\|_\alpha^\alpha} - e^{-\|uf\|_\alpha^\alpha - \|vg\|_\alpha^\alpha}. \qquad (3.1)$$

We can find the measure of dependence $U_{\xi,\eta}$ in Section 4.7 of Samorodnitsky and Taqqu [19], as well as

$$I_{\xi,\eta}(u,v) = \|uf + vg\|_\alpha^\alpha - \|uf\|_\alpha^\alpha - \|vg\|_\alpha^\alpha \qquad (3.2)$$

and $-I_{\xi,\eta}(1,-1)$, which is called the *codifference*. Observe that

$$|U_{\xi,\eta}(u,v)| \leq |I_{\xi,\eta}(u,v)|, \qquad u, v \in \mathbb{R}, \qquad (3.3)$$

due to the inequality $|e^{-x} - e^{-y}| \leq |x-y|$, $x, y > 0$. The next result, borrowed from Delbeke and Segers [10], provides various bounds on $U_{\xi,\eta}$. For notational simplicity, we shall write $\int h \, d\mu$ for the integral $\int_S h(s) \mu(ds)$ below, whenever this is not confusing.

**Lemma 3.1.** *Let $\alpha \in (0,2)$ and the function $U_{\xi,\eta}$ be defined by (3.1). Then, for all $u, v \in \mathbb{R}$,*

$$|U_{\xi,\eta}(u,v)| \leq 2|uv|^{\alpha/2} [\xi, \eta]_2, \qquad (3.4)$$

$$|U_{\xi,\eta}(u,v)| \leq 2|uv|^{\alpha/2} e^{-(|u|^{\alpha/2} \|\xi\|_\alpha^{\alpha/2} - |v|^{\alpha/2} \|\eta\|_\alpha^{\alpha/2})^2} [\xi, \eta]_2, \qquad (3.5)$$

$$|U_{\xi,\eta}(u,v)| \leq 2|uv|^{\alpha/2} e^{-2(\|\xi\|_\alpha^{\alpha/2} \|\eta\|_\alpha^{\alpha/2} - [\xi,\eta]_2)|uv|^{\alpha/2}} [\xi, \eta]_2. \qquad (3.6)$$

**Proof.** Inequality (3.4) follows from (3.3) and

$$|I_{\xi,\eta}(u,v)| \leq 2|uv|^{\alpha/2} \int |fg|^{\alpha/2} \, d\mu = 2|uv|^{\alpha/2} [\xi, \eta]_2, \qquad (3.7)$$

which follows from relation (5.8) below. As for inequalities (3.5) and (3.6), by using $|e^x - 1| = |e^x - e^0| \leq e^{|x|}|x|$, $x \in \mathbb{R}$, and (3.7), we obtain

$$|U_{\xi,\eta}(u,v)| = e^{-\|uf\|_\alpha^\alpha - \|vg\|_\alpha^\alpha} |e^{-I_{\xi,\eta}(u,v)} - 1| \leq e^{-\|uf\|_\alpha^\alpha - \|vg\|_\alpha^\alpha} |I_{\xi,\eta}(u,v)| e^{|I_{\xi,\eta}(u,v)|}$$

$$\leq 2|uv|^{\alpha/2}[\xi,\eta]_2 e^{-\|uf\|_\alpha^\alpha - \|vg\|_\alpha^\alpha + 2|uv|^{\alpha/2} \int |fg|^{\alpha/2} d\mu}$$

$$\leq 2|uv|^{\alpha/2} e^{-(|u|^{\alpha/2}\|f\|_\alpha^{\alpha/2} - |v|^{\alpha/2}\|g\|_\alpha^{\alpha/2})^2} e^{-2(\|f\|_\alpha^{\alpha/2}\|g\|_\alpha^{\alpha/2} - \int |fg|^{\alpha/2} d\mu)|uv|^{\alpha/2}} [\xi,\eta]_2.$$

This yields both inequalities (3.5) and (3.6) since $\int |fg|^{\alpha/2} \, d\mu \leq \|f\|_\alpha^{\alpha/2} \|g\|_\alpha^{\alpha/2}$ and the exponents of the two exponentials are negative. □

The following two results concern the partial derivatives of the function $U_{\xi,\eta}$. We shall use the notation $a^{\langle p \rangle} = \text{sign}(a)|a|^p$ with $a, p \in \mathbb{R}$. The next lemma follows from Lemma 5.1 below.

**Lemma 3.2.** *Let $\alpha \in (1,2)$ and $S\alpha S$ random variables $\xi, \eta$ be given by (2.1). Then, for $u, v \in \mathbb{R}$,*

$$\frac{\partial U_{\xi,\eta}}{\partial u}(u,v) = -\alpha \int (uf + vg)^{\langle \alpha-1 \rangle} f \, d\mu \, e^{-\|uf+vg\|_\alpha^\alpha}$$

$$+ \alpha \int (uf)^{\langle \alpha-1 \rangle} f \, d\mu \, e^{-\|uf\|_\alpha^\alpha - \|vg\|_\alpha^\alpha} \tag{3.8}$$

*and*

$$\frac{\partial^2 U_{\xi,\eta}}{\partial u \, \partial v}(u,v) = -\alpha(\alpha-1) \int |uf+vg|^{\alpha-2} fg \, d\mu \, e^{-\|uf+vg\|_\alpha^\alpha}$$

$$+ \alpha^2 \int (uf+vg)^{\langle \alpha-1 \rangle} f \, d\mu \int (uf+vg)^{\langle \alpha-1 \rangle} g \, d\mu \, e^{-\|uf+vg\|_\alpha^\alpha} \tag{3.9}$$

$$- \alpha^2 \int (uf)^{\langle \alpha-1 \rangle} f \, d\mu \int (vg)^{\langle \alpha-1 \rangle} g \, d\mu \, e^{-\|uf\|_\alpha^\alpha - \|vg\|_\alpha^\alpha},$$

*provided, for the relation (3.9), that*

$$\int |uf + vg|^{\alpha-2} |f||g| \, d\mu < \infty. \tag{3.10}$$

**Remark 3.1.** When $s \in S$ is such that $f(s) = 0$ or $g(s) = 0$, we assume in (3.10) and (3.9) that $|uf(s) + vg(s)|^{\alpha-2}|f(s)||g(s)| = 0$. This expression is clearly zero if one of $f(s)$ or $g(s)$ is zero; if both are zero, then their values do not contribute in the definition (2.1) of $\xi$ and $\eta$.

In the following lemma, we provide bounds for the partial derivatives in Lemma 3.2, using Assumption (A2).

**Lemma 3.3.** *Let $\alpha \in (1,2)$ and $\xi, \eta$ be two $S\alpha S$ random variables given by (2.1). Under Assumption (A2),*

$$\left| \frac{\partial U_{\xi,\eta}}{\partial u}(u,v) \right| \leq C(|v|^{\alpha-1}[\eta,\xi]_1^* e^{-\varepsilon_2(\|u\xi\|_\alpha^\alpha + \|v\eta\|_\alpha^\alpha)} + |u|^{\alpha-1}|U_{\xi,\eta}(u,v)|), \tag{3.11}$$

$$\left|\frac{\partial^2 U_{\xi,\eta}}{\partial u\,\partial v}(u,v)\right| \leq C(U_1(u,v) + U_2(u,v) + U_3(u,v)), \tag{3.12}$$

where the constant $C$ may depend on $\alpha, \varepsilon_2, \|\xi\|_\alpha$ and $\|\eta\|_\alpha$, and

$$U_1(u,v) = \int ||uf| - |vg||^{\alpha-2} |f||g|\,\mathrm{d}\mu\, \mathrm{e}^{-\varepsilon_2(\|u\xi\|_\alpha^\alpha + \|v\eta\|_\alpha^\alpha)}, \tag{3.13}$$

$$U_2(u,v) = (|u|^{2\alpha-2} + |v|^{2\alpha-2})\mathrm{e}^{-\varepsilon_2(\|u\xi\|_\alpha^\alpha + \|v\eta\|_\alpha^\alpha)}[\xi,\eta]_1, \tag{3.14}$$

$$U_3(u,v) = |u|^{\alpha-1}|v|^{\alpha-1}|U_{\xi,\eta}(u,v)|. \tag{3.15}$$

**Proof.** We consider only the bound (3.12), which is more difficult to prove. We shall denote by $C$ a generic constant which may change from occurrence to occurrence and also by $\widetilde{U}_1, \widetilde{U}_2$ and $\widetilde{U}_3$ the three terms on the right-hand side of (3.9). Since $\alpha - 2 < 0$, by using $|uf + vg| \geq ||uf| - |vg||$ and Assumption (A2), we obtain

$$|\widetilde{U}_1| \leq C \int ||uf| - |vg||^{\alpha-2} |f||g|\,\mathrm{d}\mu\, \mathrm{e}^{-\varepsilon_2(\|uf\|_\alpha^\alpha + \|vg\|_\alpha^\alpha)} = CU_1(u,v).$$

To bound $|\widetilde{U}_2 + \widetilde{U}_3|$, add and subtract to $\widetilde{U}_2$ a similar term, where the first integral is replaced by $\int (uf)^{\langle \alpha-1 \rangle} f\,\mathrm{d}\mu$, and to $\widetilde{U}_3$ a similar term, where the exponential is replaced by $\mathrm{e}^{-\|uf+vg\|_\alpha^\alpha}$. Then, by twice applying the triangle inequality,

$$|\widetilde{U}_2 + \widetilde{U}_3| \leq C \int |u|^{\alpha-1}|f|^\alpha\,\mathrm{d}\mu \int |v|^{\alpha-1}|g|^\alpha\,\mathrm{d}\mu |U_{\xi,\eta}(u,v)|$$

$$+ C\int |uf+vg|^{\alpha-1}|g|\,\mathrm{d}\mu\, \mathrm{e}^{-\|uf+vg\|_\alpha^\alpha} \int |(uf+vg)^{\langle\alpha-1\rangle} - (uf)^{\langle\alpha-1\rangle}||f|\,\mathrm{d}\mu$$

$$+ C\int |u|^{\alpha-1}|f|^\alpha\,\mathrm{d}\mu\, \mathrm{e}^{-\|uf+vg\|_\alpha^\alpha} \int |(uf+vg)^{\langle\alpha-1\rangle} - (vg)^{\langle\alpha-1\rangle}||g|\,\mathrm{d}\mu.$$

Now, by using inequality (5.7) below and Hölder's inequality, we obtain

$$|\widetilde{U}_2 + \widetilde{U}_3| \leq C \|f\|_\alpha^\alpha \|g\|_\alpha^\alpha |u|^{\alpha-1}|v|^{\alpha-1}|U_{\xi,\eta}(u,v)|$$

$$+ C(|u|^{\alpha-1}\|f\|_\alpha^{\alpha-1}\|g\|_\alpha + |v|^{\alpha-1}\|g\|_\alpha^\alpha)\mathrm{e}^{-\|uf+vg\|_\alpha^\alpha} |v|^{\alpha-1}\int |g|^{\alpha-1}|f|\,\mathrm{d}\mu$$

$$+ C|u|^{\alpha-1}\|f\|_\alpha^\alpha \mathrm{e}^{-\|uf+vg\|_\alpha^\alpha}|u|^{\alpha-1}\int |f|^{\alpha-1}|g|\,\mathrm{d}\mu$$

$$\leq C'(U_3(u,v) + U_2(u,v)),$$

using Assumption (A2) and the relation $|uv|^{\alpha-1} \leq C(|u|^{2\alpha-2} + |v|^{2\alpha-2})$. □

In the following three results, we give bounds on various integrals of $U_{\xi,\eta}$ and integrals of its partial derivatives.

**Lemma 3.4.** *Let $\alpha \in (0,2)$ and $\xi, \eta$ be $S\alpha S$ random variables that satisfy Assumption* (A1). *Then*

$$\int_{\mathbb{R}} \int_{\mathbb{R}} |U_{\xi,\eta}(u,v)| \, du \, dv \leq C[\xi,\eta]_2, \tag{3.16}$$

*where constant $C$ depends on $\varepsilon_1, \alpha, \|\xi\|_\alpha$ and $\|\eta\|_\alpha$.*

**Proof.** We consider the integral over $(0,\infty) \times (0,\infty)$ only and examine it over four regions. Over $(0,1) \times (0,1)$, by using (3.4), we have

$$\int_0^1 \int_0^1 |U_{\xi,\eta}(u,v)| \, du \, dv \leq 2 \int_0^1 \int_0^1 u^{\alpha/2} v^{\alpha/2} \, du \, dv [\xi,\eta]_2 = C[\xi,\eta]_2.$$

Over $(1,\infty) \times (1,\infty)$, by using (3.6) and Assumption (A1), we can bound the integral by

$$2 \int_1^\infty \int_1^\infty e^{-\varepsilon_1 \|f\|_\alpha^{\alpha/2} \|g\|_\alpha^{\alpha/2} u^{\alpha/2} v^{\alpha/2}} u^{\alpha/2} v^{\alpha/2} \, du \, dv [\xi,\eta]_2.$$

Observe that the integral here is finite, since we can bound $e^{-\varepsilon_1 \|f\|_\alpha^{\alpha/2} \|g\|_\alpha^{\alpha/2} u^{\alpha/2} v^{\alpha/2}}$ up to a constant by $(uv)^{-p}$ for arbitrarily large $p > 0$. Over $(0,1) \times (1,\infty)$, by using (3.5), we have as a bound,

$$2 \int_0^1 \int_1^\infty e^{-(u^{\alpha/2} \|f\|_\alpha^{\alpha/2} - v^{\alpha/2} \|g\|_\alpha^{\alpha/2})^2} u^{\alpha/2} v^{\alpha/2} \, du \, dv [\xi,\eta]_2.$$

The integral here is again finite. A similar bound holds over the region $(1,\infty) \times (0,1)$ and hence the result (3.16) is valid. $\square$

**Lemma 3.5.** *Let $\alpha \in (1,2)$ and $\xi, \eta$ be $S\alpha S$ random variables that satisfy Assumptions* (A1) *and* (A2). *Also, let*

$$F(u) = \begin{cases} |u|^{-\beta}, & \text{if } |u| < 1/b, \\ b^{\beta-1} |u|^{-1}, & \text{if } |u| \geq 1/b, \end{cases} \tag{3.17}$$

*where $b \geq 1$ and $\beta \in (0, \alpha/2)$. Then,*

$$\int_{\mathbb{R}} \int_{\mathbb{R}} \left| \frac{\partial U_{\xi,\eta}}{\partial u}(u,v) \right| F(u) \, du \, dv \leq C b^{\beta-1} \log(b+1) ([\eta,\xi]_1^* + [\xi,\eta]_2), \tag{3.18}$$

$$\int_{\mathbb{R}} \int_{\mathbb{R}} \left| \frac{\partial^2 U_{\xi,\eta}}{\partial u \, \partial v}(u,v) \right| F(u) F(v) \, du \, dv \leq C b^{2\beta-\alpha} ([\xi,\eta]_1 + [\xi,\eta]_2), \tag{3.19}$$

*where the constant $C$ depends on $\alpha, \beta, \varepsilon_1, \varepsilon_2, \|\xi\|_\alpha$ and $\|\eta\|_\alpha$. In particular, the derivative $\frac{\partial^2}{\partial v \, \partial u} U_{\xi,\eta}$ is well defined a.e. $du \, dv$.*

**Proof.** Consider first inequality (3.19), which is more difficult to prove. By Lemma 3.3, it is enough to show that the bound (3.19) holds when $\frac{\partial^2}{\partial u \partial v} U_{\xi,\eta}$ is replaced by $U_1$, by $U_2$ and by $U_3$. Since $U_1$, $U_2$, $U_3$ and $F$ are even functions in $u$ and $v$, it is enough to consider the integral over $(0,\infty) \times (0,\infty)$. Let us denote this integral with $U_1$, $U_2$ and $U_3$ by $I_1$, $I_2$ and $I_3$, respectively. That is, $I_k = \int_0^\infty \int_0^\infty U_k(u,v) F(u) F(v) \, du \, dv$, $k = 1, 2, 3$. As in the proof of Lemma 3.3, we shall denote a generic constant by $C$.

*Bounding $I_1$*: We have

$$I_1 \leq \int_S |f(s)||g(s)| \int_0^\infty \int_0^\infty ||f(s)|u - |g(s)|v|^{\alpha-2} F(u) F(v) \, du \, dv \, \mu(ds)$$

$$=: \int_S |f(s)||g(s)| I_1(s) \mu(ds).$$

By using Lemma 5.2 below, when $f(s) \neq 0$ and $g(s) \neq 0$, we obtain

$$I_1(s) = \int_0^\infty \int_0^\infty \left| \left| \frac{f(s)}{g(s)} \right| u - v \right|^{\alpha-2} F(u) F(v) \, du \, dv |g(s)|^{\alpha-2}$$

$$\leq C b^{2\beta-\alpha} \left(1 + \left| \frac{f(s)}{g(s)} \right|^{\alpha-2}\right) |g(s)|^{\alpha-2}$$

$$= C b^{2\beta-\alpha}(|f(s)|^{\alpha-2} + |g(s)|^{\alpha-2}).$$

In view of the remark following Lemma 3.2, we conclude that $I_1 \leq C b^{2\beta-\alpha}[\xi,\eta]_1$.

*Bounding $I_2$*: We consider $I_2$ over four regions $(0,1/b) \times (0,1/b)$, $(1/b,\infty) \times (1/b,\infty)$, $(0,1/b) \times (1/b,\infty)$ and $(1/b,\infty) \times (0,1/b)$, and denote the corresponding integrals by $I_{2,1}$, $I_{2,2}$, $I_{2,3}$ and $I_{2,4}$. We shall also consider $U_2$ with only the term $|u|^{2\alpha-2}$ (in the case of $|v|^{2\alpha-2}$, the bound is obtained by symmetry). Over $(0,1/b) \times (0,1/b)$, after the changes of variables $u \to u/b, v \to v/b$, we obtain

$$I_{2,1} \leq \int_0^{1/b} \int_0^{1/b} u^{2\alpha-2} u^{-\beta} v^{-\beta} \, du \, dv [\xi,\eta]_1 \leq C b^{-2\alpha+2\beta}[\xi,\eta]_1$$

since $-\beta + 1 > 0$ and $2\alpha - 2 - \beta + 1 > 0$. Over $(1/b,\infty) \times (1/b,\infty)$,

$$I_{2,2} = \int_{1/b}^\infty \int_{1/b}^\infty u^{2\alpha-2} e^{-\varepsilon_2 |u|^\alpha \|f\|_\alpha^\alpha - \varepsilon_2 |v|^\alpha \|g\|_\alpha^\alpha} b^{\beta-1} u^{-1} b^{\beta-1} v^{-1} \, du \, dv [\xi,\eta]_1$$

$$\leq b^{2\beta-2} \int_0^\infty u^{2\alpha-3} e^{-\varepsilon_2 |u|^\alpha \|f\|_\alpha^\alpha} \, du \int_{1/b}^\infty e^{-\varepsilon_2 |v|^\alpha \|g\|_\alpha^\alpha} v^{-1} \, dv [\xi,\eta]_1$$

$$\leq C b^{2\beta-2} \log(b+1)[\xi,\eta]_1$$

since $2\alpha - 2 > 0$. Over $(0,1/b) \times (1/b,\infty)$,

$$I_{2,3} \leq \int_0^{1/b} \int_{1/b}^\infty u^{2\alpha-2} e^{-\varepsilon_2 |v|^\alpha \|g\|_\alpha^\alpha} u^{-\beta} b^{\beta-1} v^{-1} \, du \, dv [\xi,\eta]_1$$

$$= b^{-2\alpha+2+\beta-1}b^{\beta-1}\int_0^1 u^{2\alpha-2-\beta}\,du \int_{1/b}^\infty e^{-\varepsilon_2|v|^\alpha\|g\|_\alpha^\alpha}v^{-1}\,dv[\xi,\eta]_1$$

$$\leq Cb^{-2\alpha+2\beta}\log(b+1)[\xi,\eta]_1.$$

Similarly, over $(1/b,\infty)\times(0,1/b)$,

$$I_{2,4}\leq \int_{1/b}^\infty \int_0^{1/b} u^{2\alpha-2}e^{-\varepsilon_2|u|^\alpha\|f\|_\alpha^\alpha}b^{\beta-1}u^{-1}v^{-\beta}\,du\,dv[\xi,\eta]_1$$

$$\leq \int_0^\infty u^{2\alpha-3}e^{-\varepsilon_2|u|^\alpha\|f\|_\alpha^\alpha}\,du \int_0^{1/b} v^{-\beta}\,dv\,b^{\beta-1}[\xi,\eta]_1 = Cb^{2\beta-2}[\xi,\eta]_1.$$

We conclude that $I_2 \leq Cb^{2\beta-2}\log(b+1)[\xi,\eta]_1$.

*Bounding $I_3$:* We shall here use inequalities (3.4)–(3.6). As in the case of $I_2$, denote the integral $I_3$ over the same four regions by $I_{3,1}$, $I_{3,2}$, $I_{3,3}$ and $I_{3,4}$, respectively. Then, by using (3.4), we have

$$I_{3,1}\leq C\int_0^{1/b}\int_0^{1/b} u^{\alpha-1}v^{\alpha-1}u^{\alpha/2}v^{\alpha/2}u^{-\beta}v^{-\beta}\,du\,dv[\xi,\eta]_2 \leq Cb^{-3\alpha+2\beta}[\xi,\eta]_2.$$

By using (3.6) and Assumption (A1), and making the changes of variables $u = b^{-1}x^{2/\alpha}, v = b^{-1}y^{2/\alpha}$ below, we obtain

$$I_{3,2}\leq Cb^{2\beta-2}\int_{1/b}^\infty\int_{1/b}^\infty u^{\alpha-1}v^{\alpha-1}u^{\alpha/2}v^{\alpha/2}e^{-2\varepsilon_1\|f\|_\alpha^{\alpha/2}\|g\|_\alpha^{\alpha/2}|uv|^{\alpha/2}}u^{-1}v^{-1}\,du\,dv[\xi,\eta]_2$$

$$\leq Cb^{2\beta-3\alpha}\int_1^\infty\int_1^\infty e^{-2\varepsilon_1\|f\|_\alpha^{\alpha/2}\|g\|_\alpha^{\alpha/2}b^{-\alpha}xy}(xy)^{2(\alpha-1)/\alpha}\,dx\,dy[\xi,\eta]_2.$$

By a further changes of variables $x = b^\alpha w$ and $wy = z$, we have

$$I_{3,2}\leq Cb^{2\beta-3\alpha+\alpha+2(\alpha-1)}\int_{b^{-\alpha}}^\infty \frac{dw}{w}\int_w^\infty e^{-2\varepsilon_1\|f\|_\alpha^{\alpha/2}\|g\|_\alpha^{\alpha/2}z}z^{2(\alpha-1)/\alpha}\,dz[\xi,\eta]_2$$

$$\leq Cb^{2\beta-2}\log(b+1)[\xi,\eta]_2,$$

by splitting and bounding $\int_{b^{-\alpha}}^\infty\int_w^\infty \leq \int_{b^{-\alpha}}^1\int_0^\infty + \int_1^\infty\int_w^\infty$. Turning to $I_{3,3}$, this time using (3.5) and making the changes of variables $u = x^{2/\alpha}b^{-1}, v = y^{2/\alpha}b^{-1}$ below, we obtain

$$I_{3,3}\leq Cb^{\beta-1}\int_0^{1/b}\int_{1/b}^\infty u^{\alpha-1}v^{\alpha-1}u^{\alpha/2}v^{\alpha/2}e^{-(u^{\alpha/2}\|f\|_\alpha^{\alpha/2}-v^{\alpha/2}\|g\|_\alpha^{\alpha/2})^2}u^{-\beta}v^{-1}\,du\,dv[\xi,\eta]_2$$

$$\leq Cb^{2\beta-3\alpha}\int_0^1\int_1^\infty x^{2(\alpha-\beta)/\alpha}y^{2(\alpha-1)/\alpha}e^{-(xb^{-\alpha/2}\|f\|_\alpha^{\alpha/2}-yb^{-\alpha/2}\|g\|_\alpha^{\alpha/2})^2}\,dx\,dy[\xi,\eta]_2.$$

By making another change of variables $yb^{-\alpha/2} = w$, it is easy to see that the last expression can be bounded (up to a constant) by $b^{2\beta-3\alpha}b^{(\alpha/2)+\alpha-1}[\xi,\eta]_2 = b^{2\beta-((3\alpha)/2)-1}[\xi,\eta]_2$.

One obtains the same bound for $I_{3,4}$ by symmetry. We can now conclude that $I_3 \leq Cb^{2\beta-2}\log(b+1)[\xi,\eta]_2$.

The result (3.19) of the lemma follows from the bounds obtained for $I_1$, $I_2$ and $I_3$. The inequality (3.18) can be shown in a similar way by using (3.11). □

**Lemma 3.6.** *Let $\alpha \in (1,2)$ and $\xi, \eta$ be $S\alpha S$ random variables (2.1) satisfying Assumptions (A1) and (A2). Also, let*

$$G(u) = \begin{cases} |u|^{-\beta}, & \text{if } |u| < 1, \\ 1, & \text{if } |u| \geq 1, \end{cases} \tag{3.20}$$

*where $\beta \in (0, \alpha/2)$. Then,*

$$\int_{\mathbb{R}} \int_{\mathbb{R}} \left| \frac{\partial U_{\xi,\eta}}{\partial u}(u,v) \right| G(u) \, du \, dv \leq C < \infty, \tag{3.21}$$

$$\int_{\mathbb{R}} \int_{\mathbb{R}} \left| \frac{\partial^2 U_{\xi,\eta}}{\partial u \, \partial v}(u,v) \right| G(u) G(v) \, du \, dv \leq C < \infty, \tag{3.22}$$

*where the constant $C$ depends on $\alpha, \beta, \varepsilon_1, \varepsilon_2, \|\xi\|_\alpha$ and $\|\eta\|_\alpha$.*

**Proof.** We prove only the bound (3.22). The generic constants $C$ and $C'$ below are as in the statement of the lemma. As in the proof of Lemma 3.5, we consider the integral (3.22) over $(0,\infty) \times (0,\infty)$ only and denote it by $J_1$, $J_2$ and $J_3$ when $\frac{\partial^2}{\partial u \, \partial v} U_{\xi,\eta}$ is replaced by $U_1$, $U_2$ and $U_3$ in Lemma 3.3, respectively. To bound $J_1$, we can apply the bound for the integral $I_1$ in the proof of Lemma 3.5 with $b=1$ because the difference between the functions $F$ and $G$ can be accounted for by bounding the term $e^{-\varepsilon_2(\|uf\|_\alpha^\alpha + \|vg\|_\alpha^\alpha)}$ in $U_1$ by $C|u|^{-1}|v|^{-1}$ when $|u| > 1$ or $|v| > 1$. We therefore obtain $J_1 \leq C[\xi,\eta]_1 \leq C(\|\xi\|_\alpha^{\alpha-1}\|\eta\|_\alpha + \|\eta\|_\alpha^{\alpha-1}\|\xi\|_\alpha) \leq C' < \infty$. In the case of $J_2$, we have

$$J_2 = \int_0^\infty \int_0^\infty (|u|^{2\alpha-2} + |v|^{2\alpha-2}) e^{-\varepsilon_2(\|uf\|_\alpha^\alpha + \|vg\|_\alpha^\alpha)} G(u) G(v) \, du \, dv \, [\xi,\eta]_1$$

and by bounding the exponential by $Ce^{-\varepsilon_3(\|uf\|_\alpha^\alpha + \|vg\|_\alpha^\alpha)} |u|^{-1}|v|^{-1}$ when $|u| > 1$ or $|v| > 1$, we can again account for the difference between the functions $F$ and $G$. We can therefore again apply the bound for the integral $I_2$ in the proof of Lemma 3.3 with $b+1$ and obtain $J_2 \leq C[\xi,\eta]_1 \leq C' < \infty$, as above. The integral

$$J_3 = \int_0^\infty \int_0^\infty |u|^{\alpha-1}|v|^{\alpha-1}|U_{\xi,\eta}(u,v)| G(u) G(v) \, du \, dv$$

can be bounded by using techniques from the proof of Lemma 3.4, together with the bounds (3.4)–(3.6) for $U_{\xi,\eta}$ and the inequality $[\xi,\eta]_2 \leq \|\xi\|_\alpha^{\alpha/2}\|\eta\|_\alpha^{\alpha/2} = C < \infty$. □

## 4. Proofs of the main results

We shall prove here Propositions 2.1 and 2.2 and Theorems 2.1–2.4.

**Proof of Proposition 2.1.** The vector $(\xi, \eta)$ in (2.1) has the so-called minimal representation

$$(\xi, \eta) \stackrel{d}{=} \left( \int_{\widetilde{S}} \widetilde{f}(\widetilde{s}) \widetilde{M}(d\widetilde{s}), \int_{\widetilde{S}} \widetilde{g}(\widetilde{s}) \widetilde{M}(d\widetilde{s}) \right), \tag{4.1}$$

where $\widetilde{M}$ is an $S\alpha S$ random measure on $\widetilde{S}$ with the control measure $\widetilde{\mu}$. For more information on minimal representations, see Rosiński [18] or [17], Section 2. It is therefore enough to show that the quantities $[\xi, \eta]_1^*$, $[\xi, \eta]_1$ and $[\xi, \eta]_2$ in (2.4)–(2.5) remain the same when either the representation (2.1) is used or when the minimal representation (4.1) is used. We consider only the quantity $[\xi, \eta]_2$ in (2.5).

We may assume, without loss of generality, that $\{s: f(s) \neq 0 \text{ or } g(s) \neq 0\} = S$ $\mu$-a.e. (otherwise, the set $\{s: f(s) = 0, g(s) = 0\} = S$ can be eliminated from the representation (2.1) without changing $[\xi, \eta]_2$). By Remark 2.5 in Rosiński [17], there are maps $\phi: S \mapsto \widetilde{S}$ and $h: S \mapsto \mathbb{R} \setminus \{0\}$ such that

$$f(s) = h(s)\widetilde{f}(\phi(s)), \qquad g(s) = h(s)\widetilde{g}(\phi(s)) \tag{4.2}$$

a.e. $\mu(ds)$ and

$$\widetilde{\mu} = \mu_h \circ \phi^{-1}, \tag{4.3}$$

where $\mu_h(ds) = |h(s)|^\alpha \mu(ds)$. By using (4.2) and (4.3), and making a change of variables, we have

$$\int_S |f(s)g(s)|^{\alpha/2} \mu(ds) = \int_S |\widetilde{f}(\phi(s))\widetilde{g}(\phi(s))|^{\alpha/2} \mu_h(ds) = \int_{\widetilde{S}} |\widetilde{f}(\widetilde{s})\widetilde{g}(\widetilde{s})|^{\alpha/2} \widetilde{\mu}(d\widetilde{s}).$$

This proves the result for the quantity $[\xi, \eta]_2$ in (2.5). □

**Proof of Proposition 2.2.** Let $(T_n a)(x) = a(n-x)$, $n \in \mathbb{Z}, x \in \mathbb{R}$. In view of (2.3), $\xi_0$ and $\xi_n$ have kernels $f = T_0 a$ and $g = T_n a$. That $\xi_0$ and $\xi_n$ satisfy Assumption (A1) follows from the facts that $\|T_n a\|_\alpha = \|T_0 a\|_\alpha$ and $\int_\mathbb{R} |(T_0 a)(x)(T_n a)(x)|^{\alpha/2} dx \to 0$ as $n \to \infty$. To see why the last integral converges to 0, first write it as

$$\int_\mathbb{R} |a(-x)a(n-x)|^{\alpha/2} dx = \int_l^\infty |a(-x)a(n-x)|^{\alpha/2} dx + \int_{-\infty}^l |a(-x)a(n-x)|^{\alpha/2} dx$$

$$\leq \|a(-\cdot)1_{(l,\infty)}\|_\alpha^{\alpha/2} \|a(-\cdot)\|_\alpha^{\alpha/2} + \|a(-\cdot)\|_\alpha^{\alpha/2} \|a(\cdot)1_{(n-l,\infty)}\|_\alpha^{\alpha/2},$$

by the Cauchy–Schwarz inequality and a change of variables. Now, observe that the first term in the bound is arbitrarily small for large fixed $l$ and that the second term converges to 0 for fixed $l$ as $n \to \infty$.

We now show that they also satisfy Assumption (A2). We consider only the case $\alpha \in [1,2)$; the case $\alpha \in (0,1)$ can be proven by considering $\|\cdot\|_\alpha^\alpha$ below. For some $l \in \mathbb{R}$, we have, by $2(x+y)^p \geq x^p + y^p$, $x,y,p > 0$, and by Minkowski's inequality applied to the norm $\|\cdot\|_\alpha$,

$$\begin{aligned}
\|uT_0a + vT_na\|_\alpha &= \|(uT_0a + vT_na)1_{(-\infty,l)} + (uT_0a + vT_na)1_{(l,\infty)}\|_\alpha \\
&\geq \tfrac{1}{2}\|(uT_0a + vT_na)1_{(-\infty,l)}\|_\alpha + \tfrac{1}{2}\|(uT_0a + vT_na)1_{(l,\infty)}\|_\alpha \\
&\geq \tfrac{1}{2}(\|uT_0a1_{(-\infty,l)}\|_\alpha - \|vT_na1_{(-\infty,l)}\|_\alpha) \\
&\quad + \tfrac{1}{2}(\|vT_na1_{(l,\infty)}\|_\alpha - \|uT_0a1_{(l,\infty)}\|_\alpha) \\
&= \tfrac{1}{2}(\|T_0a1_{(-\infty,l)}\|_\alpha - \|T_0a1_{(l,\infty)}\|_\alpha)|u| \\
&\quad + \tfrac{1}{2}(\|T_0a1_{(l-n,\infty)}\|_\alpha - \|T_0a1_{(-\infty,l-n)}\|_\alpha)|v|.
\end{aligned} \quad (4.4)$$

Since $\|T_0a\|_\alpha = \lim_{m\to\infty} \|T_0a1_{(-\infty,m)}\|_\alpha$, there is a constant $M > 0$ such that, for all $m \geq M$,

$$\|T_0a1_{(-\infty,m)}\|_\alpha - \|T_0a1_{(m,\infty)}\|_\alpha \geq \tfrac{1}{2}\|T_0a\|_\alpha,$$
$$\|T_0a1_{(-m,\infty)}\|_\alpha - \|T_0a1_{(-\infty,-m)}\|_\alpha \geq \tfrac{1}{2}\|T_0a\|_\alpha.$$

By applying these inequalities to (4.4), by fixing $l \geq M$ and taking $n$ such that $n - l \geq M$ or $n \geq (M+l)$, we obtain

$$\|uT_0a + vT_na\|_\alpha \geq \frac{1}{4}(\|uT_0a\|_\alpha + \|vT_na\|_\alpha) \geq \frac{1}{4\cdot 2^{1/\alpha}}(\|uT_0a\|_\alpha^\alpha + \|vT_na\|_\alpha^\alpha)^{1/\alpha},$$

which is Assumption (A2). □

We shall denote the Fourier transform and its inverse by

$$\begin{aligned}
(\mathcal{F}h)(u) &= \int_{\mathbb{R}^n} e^{iu\cdot x} h(x)\, dx, \\
(\mathcal{F}^{-1}k)(x) &= (\mathcal{G}k)(x) = \frac{1}{(2\pi)^n} \int_{\mathbb{R}^n} e^{-ix\cdot u} k(u)\, du
\end{aligned} \quad (4.5)$$

and we shall also use the measure of dependence $U_{\xi,\eta}$, which was analyzed in Section 3. The following two lemmas are used in the proofs of Theorems 2.1–2.4.

**Lemma 4.1.** *Let $K$ be a function satisfying Assumptions (K1) and (K2) with $x_0 = b \geq 1$. Set $K_b(x) = K(x)1_{\{|x|>b\}}$ and $K_{b,n}(x) = K(x)1_{\{b<|x|<n\}}$, $x \in \mathbb{R}, n \in \mathbb{N}$. Then,*

$$\left|\left(\mathcal{G}\frac{K_b(x)}{x}\right)(u)\right| \leq C \begin{cases} |u|^{-\beta}, & \text{if } |u| < 1/b, \\ b^{\beta-1}|u|^{-1}, & \text{if } u \in \mathbb{R}, \end{cases} \quad (4.6)$$

where the constant $C$ does not depend on $b$. Moreover, for fixed $b$,

$$\left|\left(\mathcal{G}\frac{K_{b,n}(x)}{x}\right)(u)\right| \leq C \begin{cases} |u|^{-\beta}, & \text{if } |u| < 1/b, \\ |u|^{-1}, & \text{if } u \in \mathbb{R}, \end{cases} \quad (4.7)$$

where the constant $C$ does not depend on $n$.

**Proof.** The second inequality in (4.6) follows from

$$2\pi\left|\left(\mathcal{G}\frac{K_b(x)}{x}\right)(u)\right| \leq \left|\left(\int_{-\infty}^{-b} + \int_{b}^{\infty}\right)\frac{K(x)}{x}e^{-ixu}\,dx\right|$$
$$\leq \frac{4(|K(b)| + |K(-b)|)}{b|u|} \leq Cb^{\beta-1}|u|^{-1},$$

where we used Assumptions (K1) and (K2), together with the inequality

$$\left|\int_c^d e^{ixu} f(x)\,dx\right| \leq \frac{4}{|u|} \max_{x \in [c,d]} |f(x)|, \quad (4.8)$$

valid for $c < d$, $u \in \mathbb{R}$ and a function $f$ monotonic on $[c,d]$. (The inequality (4.8) is a consequence of the second mean value theorem; see, e.g., in Kawata [12], page 24). As for the first inequality in (4.6), when $0 < |u| < 1/b$, we obtain

$$2\pi\left|\left(\mathcal{G}\frac{K_b(x)}{x}\right)(u)\right| \leq \int_\mathbb{R} 1_{\{b < |x| < 1/|u|\}}\frac{|K(x)|}{|x|}\,dx + \left|\left(\int_{-\infty}^{-1} + \int_1^\infty\right)\frac{K(z/|u|)}{z}e^{-iz\,\text{sign}(u)}\,dz\right|.$$

Using the fact that $\frac{K(z/|u|)}{z/|u|}$ is monotone for $|z|/|u| > b$, hence for $|z| > 1$, and applying (4.8), we obtain

$$2\pi\left|\left(\mathcal{G}\frac{K_b(x)}{x}\right)(u)\right| \leq 2\int_b^{1/|u|} x^{\beta-1}\,dx + 4\left(K\left(\frac{1}{|u|}\right) + K\left(-\frac{1}{|u|}\right)\right) \leq C|u|^{-\beta}.$$

The bound (4.7) can be proven in a similar way. $\square$

**Lemma 4.2.** *Let $\alpha \in (1,2)$ and $\xi, \eta$ be two $S\alpha S$ random variables (2.1) that satisfy Assumptions (A1) and (A2). Then, for $(x,y) \in \mathbb{R}^2$, $x, y \neq 0$,*

$$(\mathcal{G}U_{\xi,\eta})(x,y) = \frac{1}{ix}\left(\mathcal{G}\frac{\partial}{\partial u}U_{\xi,\eta}\right)(x,y), \quad (4.9)$$

$$(\mathcal{G}U_{\xi,\eta})(x,y) = -\frac{1}{xy}\left(\mathcal{G}\frac{\partial^2}{\partial u\,\partial v}U_{\xi,\eta}\right)(x,y). \quad (4.10)$$

**Proof.** Since, by Lemma 3.6, $\frac{\partial}{\partial u}U_{\xi,\eta}, \frac{\partial^2}{\partial u\,\partial v}U_{\xi,\eta} \in L^1(\mathbb{R}^2)$ and also, by Assumption (A2) and the bound (3.11), $\lim_{u\to\infty} U_{\xi,\eta}(u,v) = 0$ and $\lim_{v\to\infty}(\frac{\partial}{\partial u}U_{\xi,\eta})(u,v) = 0$, we obtain, using the integration by parts formula, that

$$(2\pi)^2(\mathcal{G}U_{\xi,\eta})(x,y) = \int_{\mathbb{R}^2} e^{-ixu-iyv} U_{\xi,\eta}(u,v)\, du\, dv$$

$$= -\frac{1}{ix}\int_{\mathbb{R}} dv\, e^{-iyv}\int_{\mathbb{R}} d(e^{-ixu})U_{\xi,\eta}(u,v)$$

$$= \frac{1}{ix}\int_{\mathbb{R}} dv\, e^{-iyv}\int_{\mathbb{R}} du\, e^{-ixu}\frac{\partial U_{\xi,\eta}}{\partial u}(u,v) = \frac{(2\pi)^2}{ix}\left(\mathcal{G}\frac{\partial}{\partial u}U_{\xi,\eta}\right)(x,y)$$

$$= -\frac{1}{(ix)(iy)}\int_{\mathbb{R}} du\, e^{-ixu}\int_{\mathbb{R}} d(e^{-iyv})\frac{\partial U_{\xi,\eta}}{\partial u}(u,v)$$

$$= -\frac{(2\pi)^2}{xy}\left(\mathcal{G}\frac{\partial^2}{\partial u\, \partial v}U_{\xi,\eta}\right)(x,y). \qquad \square$$

We next provide the proofs of Theorems 2.1–2.4.

**Proof of Theorem 2.1.** Let $\psi_{\xi,\eta}$, $\psi_\xi$ and $\psi_\eta$ be the density functions of $(\xi,\eta)$, $\xi$ and $\eta$, respectively. Also, let $\phi_{\xi,\eta}$, $\phi_\xi$ and $\phi_\eta$ be the characteristic functions of $(\xi,\eta)$, $\xi$ and $\eta$, respectively. By using the inversion formula for the densities of $(\xi,\eta)$, $\xi$ and $\eta$, and $K, L \in L^1(\mathbb{R})$ and $U_{\xi,\eta} \in L^1(\mathbb{R}^2)$, we obtain that

$$\mathrm{Cov}(K(\xi), L(\eta)) = EK(\xi)L(\eta) - EK(\xi)EL(\eta)$$

$$= \int_{\mathbb{R}^2} K(x)L(y)(\psi_{\xi,\eta}(u,v) - \psi_\xi(u)\psi_\eta(v))\, dx\, dy$$

$$= \int_{\mathbb{R}^2} K(x)L(y)((\mathcal{G}\phi_{\xi,\eta})(x,y) - (\mathcal{G}\phi_\xi)(x)(\mathcal{G}\phi_\eta)(y))\, dx\, dy$$

$$= \int_{\mathbb{R}^2} K(x)L(y)(\mathcal{G}U_{\xi,\eta})(x,y)\, dx\, dy$$

$$= \int_{\mathbb{R}^2} (\mathcal{G}K)(u)(\mathcal{G}L)(v)U_{\xi,\eta}(u,v)\, du\, dv,$$

by Fubini's theorem. The bound (2.9) follows by using $|(\mathcal{G}K)(u)| \leq C\|K\|_1$, $|(\mathcal{G}L)(u)| \leq C\|L\|_1$ and Lemma 3.4. $\qquad \square$

**Proof of Theorem 2.2.** For $n \geq 1$, set $K_{b,n}(x) = K_b(x)1_{\{|x|\leq n\}} = K(x)1_{\{b<|x|\leq n\}}$ and $L_{b,n}(x) = L(x)1_{\{b<|x|\leq n\}}$, $x \in \mathbb{R}$. Since, by Lemma 3.6, $\frac{\partial^2}{\partial u\, \partial v}U_{\xi,\eta} \in L^1(\mathbb{R}^2)$ and $K_{b,n}, L_{b,n} \in L^1(\mathbb{R})$, we obtain, using (4.10), that

$$\mathrm{Cov}(K_{b,n}(\xi), L_{b,n}(\eta)) = \int_{\mathbb{R}^2} K_{b,n}(x)L_{b,n}(y)(\mathcal{G}U_{\xi,\eta})(x,y)\, dx\, dy$$

$$= -\int_{\mathbb{R}^2} \frac{K_{b,n}(x)}{x}\frac{L_{b,n}(y)}{y}\left(\mathcal{G}\frac{\partial^2}{\partial u\, \partial v}U_{\xi,\eta}\right)(x,y)\, dx\, dy \qquad (4.11)$$

$$= -\int_{\mathbb{R}^2}\left(\mathcal{G}\frac{K_{b,n}(x)}{x}\right)(u)\left(\mathcal{G}\frac{L_{b,n}(y)}{y}\right)(v)\frac{\partial^2 U_{\xi,\eta}}{\partial u\,\partial v}(u,v)\,\mathrm{d}u\,\mathrm{d}v.$$

Now, observe that, by using Assumption (K2) and the inequality (4.8) above, we have

$$\begin{aligned}\left(\mathcal{G}\frac{K_{b,n}(x)}{x}\right)(u) &= \frac{1}{2\pi}\int_{\mathbb{R}} e^{-\mathrm{i}ux}\frac{K(x)}{x}1_{\{b<|x|\le n\}}\,\mathrm{d}x \\ &\longrightarrow \frac{1}{2\pi}\int_{\mathbb{R}} e^{-\mathrm{i}ux}\frac{K(x)}{x}1_{\{b<|x|\}}\,\mathrm{d}x =: \left(\mathcal{G}\frac{K_b(x)}{x}\right)(u)\end{aligned} \qquad (4.12)$$

as $n\to\infty$, for all $u\in\mathbb{R}$ ($u\neq 0$). Then, by letting $n\to\infty$ in (4.11) and using (4.7) in Lemma 4.1, Lemma 3.6 and the dominated convergence theorem, we obtain that

$$\mathrm{Cov}(K_b(\xi), L_b(\eta)) = -\int_{\mathbb{R}^2}\left(\mathcal{G}\frac{K_b(x)}{x}\right)(u)\left(\mathcal{G}\frac{L_b(y)}{y}\right)(v)\frac{\partial^2 U_{\xi,\eta}}{\partial u\,\partial v}(u,v)\,\mathrm{d}u\,\mathrm{d}v. \qquad (4.13)$$

The bound (2.10) follows from (4.6) in Lemma 4.1 and (3.19) in Lemma 3.5. □

*Remark 4.1.* Note that the proof of Theorem 2.1 fails for Theorem 2.2 because it may happen that $\mathcal{G}K_b$ is not defined (this is the case, e.g., when $K_b(x) = \log_2|x|1_{\{|x|>b\}}$). This is why we consider $\mathcal{G}K_b(x)/x$ instead in (4.12). Note, also, that $(\mathcal{G}K_b(x)/x)(u)$ in (4.12) is an improper integral because we could have $K_b(x)/x \notin L^1(\mathbb{R})\cup L^2(\mathbb{R})$ and thus the usual $L^1$- or $L^2$-Fourier transforms for $K_b(x)/x$ may not be defined. What helped was the monotonicity of the function $K_b(x)/x$.

*Remark 4.2.* In order to avoid dealing with complications involving $\frac{\partial^2}{\partial u\,\partial v}U_{\xi,\eta}$, one may think of transferring (by the integration by parts formula) the derivative $\frac{\partial^2}{\partial u\,\partial v}$ from $U_{\xi,\eta}$ onto $(\mathcal{G}\cdot)$ to obtain

$$\mathrm{Cov}(K_b(\xi), L_b(\eta)) = -\int_{\mathbb{R}^2}\frac{\partial}{\partial u}\left(\mathcal{G}\frac{K_b(x)}{x}\right)(u)\frac{\partial}{\partial v}\left(\mathcal{G}\frac{L_b(y)}{y}\right)(v)U_{\xi,\eta}(u,v)\,\mathrm{d}u\,\mathrm{d}v. \qquad (4.14)$$

The problem with (4.14), however, is that the derivatives of $(\mathcal{G}\cdot)$ are not easy to manipulate. Consider, for example, the functions $K_b(x) = L_b(x) = |x|^\beta 1_{\{|x|>b\}}$ with $\beta \in (0,\alpha/2)$ which appear in the wavelet applications. We can verify that

$$\frac{\partial}{\partial u}\left(\mathcal{G}\frac{K_b(x)}{x}\right)(u) = \frac{\mathrm{i}}{\pi}\mathrm{sign}(u)\left(\beta u^{-\beta-1}\int_{bu}^\infty z^{\beta-1}\sin z\,\mathrm{d}z + b^\beta\frac{\sin bu}{u}\right)$$

and hence the right-hand side of (4.14) involves the integral with $\sin bu\sin bv$ as part of its integrand. It turns out that we cannot estimate this integral by putting absolute values on the integrand (the sign of $\sin bu\sin bv$ is important). We therefore worked with the formula (4.13) instead.

**Remark 4.3.** Although the formula (4.14) is not used here, it may be useful in other contexts. For example, if we take $\ln|\cdot|$ for the functions $K_b$ and $L_b$, we obtain

$$\operatorname{Cov}(\ln|\xi|, \ln|\eta|) = \frac{1}{4} \int_{\mathbb{R}} \int_{\mathbb{R}} U_{\xi,\eta}(u,v) \frac{\mathrm{d}u\,\mathrm{d}v}{|u||v|}. \tag{4.15}$$

The formula (4.15) is proved rigorously in Delbeke and Segers [10], where it is used to bound the covariance (4.15).

**Proof of Theorem 2.3.** For $n \geq 1$, set $K_{b,n}(x) = K(x)1_{\{b<|x|\leq n\}}$, $x \in \mathbb{R}$. Since, by Lemma 3.6, $\frac{\partial}{\partial u}U_{\xi,\eta} \in L^1(\mathbb{R}^2)$ and $K_{b,n} \in L^1(\mathbb{R})$, we obtain, using Lemma 4.2, that

$$\begin{aligned}
\operatorname{Cov}(K_{b,n}(\xi), L(\eta)) &= \int_{\mathbb{R}^2} K_{b,n}(x) L(y)(\mathcal{G}U_{\xi,\eta})(x,y)\,\mathrm{d}x\,\mathrm{d}y \\
&= \frac{1}{\mathrm{i}} \int_{\mathbb{R}^2} \frac{K_{b,n}(x)}{x} L(y) \left(\mathcal{G}\frac{\partial}{\partial u}U_{\xi,\eta}\right)(x,y)\,\mathrm{d}x\,\mathrm{d}y \\
&= \frac{1}{\mathrm{i}} \int_{\mathbb{R}^2} \left(\mathcal{G}\frac{K_{b,n}(x)}{x}\right)(u)(\mathcal{G}L)(v)\frac{\partial U_{\xi,\eta}}{\partial u}(u,v)\,\mathrm{d}u\,\mathrm{d}v.
\end{aligned} \tag{4.16}$$

By letting $n \to \infty$ in (4.16) and using (4.12), (4.7), Lemma 3.6 and the dominated convergence theorem, we obtain that

$$\operatorname{Cov}(K_b(\xi), L(\eta)) = \frac{1}{\mathrm{i}} \int_{\mathbb{R}^2} \left(\mathcal{G}\frac{K_b(x)}{x}\right)(u)(\mathcal{G}L)(v)\frac{\partial U_{\xi,\eta}}{\partial u}(u,v)\,\mathrm{d}u\,\mathrm{d}v. \tag{4.17}$$

The bound (2.11) follows from (4.6) in Lemma 4.1 and (3.18) in Lemma 3.5. □

**Proof of Theorem 2.4.** Let $x_{0,K}$ and $x_{0,L}$ be the points in Assumptions (K1)–(K3) corresponding to the functions $K$ and $L$, respectively. For $b > \max\{x_{0,K}, x_{0,L}\}$ and $b \geq 1$, set $K_b(x) = K(x)1_{\{|x|>b\}}$ and $L_b(x) = L(x)1_{\{|x|>b\}}$, $x \in \mathbb{R}$. By Assumption (K3), we have $K = K_b + \widetilde{K}$ and $L = L_b + \widetilde{L}$, with $\widetilde{K}, \widetilde{L} \in L^1(\mathbb{R}, \mathrm{d}x)$. The bound (2.12) follows by first writing

$$\begin{aligned}
|\operatorname{Cov}(K(\xi), L(\eta))| &\leq |\operatorname{Cov}(K_b(\xi), L_b(\eta))| + |\operatorname{Cov}(K_b(\xi), \widetilde{L}(\eta))| \\
&\quad + |\operatorname{Cov}(\widetilde{K}(\xi), L_b(\eta))| + |\operatorname{Cov}(\widetilde{K}(\xi), \widetilde{L}(\eta))|
\end{aligned} \tag{4.18}$$

and then using Theorems 2.1, 2.2 and 2.3 to bound the terms on the right-hand side of (4.18). □

## 5. Auxiliary results

We prove here some results which were used in Sections 2 and 3. The following lemma was used in the proof of Lemma 3.2. It involves signed powers $a^{\langle p \rangle} = |a|^p \operatorname{sign}(a)$. We write $\int h\,\mathrm{d}\mu = \int_S h(s)\mu(\mathrm{d}s)$ for notational simplicity.

**Lemma 5.1.** *Let* $f, g \in L^\alpha(S, \mu)$ *with* $\alpha \in (1, 2)$ *and set*

$$F(u,v) = \int |uf + vg|^\alpha \, \mathrm{d}\mu, \qquad u, v \in \mathbb{R}. \tag{5.1}$$

*Then, for* $u, v \in \mathbb{R}$,

$$\frac{\partial F}{\partial u}(u,v) = \alpha \int (uf + vg)^{\langle \alpha - 1 \rangle} f \, \mathrm{d}\mu,$$

$$\frac{\partial^2 F}{\partial v \, \partial u}(u,v) = \alpha(\alpha - 1) \int |uf + vg|^{\alpha - 2} fg \, \mathrm{d}\mu, \tag{5.2}$$

*provided, for the second equality in* (5.2), *that*

$$\int |uf + vg|^{\alpha - 2} |f| |g| \, \mathrm{d}\mu < \infty. \tag{5.3}$$

**Proof.** By the mean value theorem, we have, for some $u^* = u^*(u, v, h, x) \in [u, u + h]$ (we consider $h > 0$),

$$\left| \frac{1}{h}(F(u+h, v) - F(u, v)) - \alpha \int (uf + vg)^{\langle \alpha - 1 \rangle} f \, \mathrm{d}\mu \right|$$

$$= l \left| \int \frac{1}{h}(|(u+h)f + vg|^\alpha - |uf + vg|^\alpha) \, \mathrm{d}\mu - \alpha \int (uf + vg)^{\langle \alpha - 1 \rangle} f \, \mathrm{d}\mu \right|$$

$$= \alpha \left| \int ((u^*f + vg)^{\langle \alpha - 1 \rangle} - (uf + vg)^{\langle \alpha - 1 \rangle}) f \, \mathrm{d}\mu \right|$$

$$\leq 2\alpha \int |u - u^*|^{\alpha - 1} |f|^\alpha \, \mathrm{d}\mu \leq 2\alpha h^{\alpha - 1} \int |f|^\alpha \, \mathrm{d}\mu \to 0,$$

by (5.7) in Lemma 5.3, as $h \to 0$, which proves the result for the first derivative. The case of the second partial derivative of $F$ is more delicate. By the mean value theorem with $v^* = v^*(u, v, h, x) \in [v, v + h]$, we have

$$\frac{1}{h}\left( \frac{\partial F}{\partial u}(u, v+h) - \frac{\partial F}{\partial u}(u, v) \right) = \alpha \int \frac{1}{h}((uf + (v+h)g)^{\langle \alpha - 1 \rangle} - (uf + vg)^{\langle \alpha - 1 \rangle}) f \, \mathrm{d}\mu$$

$$= \alpha(\alpha - 1) \int |uf + v^*g|^{\alpha - 2} fg \, \mathrm{d}\mu.$$

Since $|uf + v^*g|^{\alpha - 2} \to |uf + vg|^{\alpha - 2}$, as $h \downarrow 0$, and, by (5.6) in Lemma 5.3 below,

$$(\alpha - 1)|uf + v^*g|^{\alpha - 2}|f||g|$$

$$= \frac{1}{h}|(uf + (v+h)g)^{\langle \alpha - 1 \rangle} - (uf + vg)^{\langle \alpha - 1 \rangle}||f| \leq 2|uf + vg|^{\alpha - 2}|f||g|,$$

the dominated convergence theorem and condition (5.3) together imply that

$$\frac{1}{h}\left(\frac{\partial F}{\partial u}(u,v+h) - \frac{\partial F}{\partial u}(u,v)\right) \xrightarrow{h\downarrow 0} \alpha(\alpha-1)\int |uf+vg|^{\alpha-2}fg\,\mathrm{d}\mu,$$

which proves the result for the second partial derivative. □

The next lemma was used in the proof of Lemma 3.5.

**Lemma 5.2.** *Let $\alpha \in (1,2)$, $r > 0$ and $F(u) = |u|^{-\beta}$ for $|u| < 1/b$ and $F(u) = b^{\beta-1}|u|^{-1}$ for $|u| \geq 1/b$, where $b \geq 1$ and $\beta \in (0, \alpha/2)$. Then,*

$$I_0(r) := \int_0^\infty \int_0^\infty |ru-v|^{\alpha-2} F(u)F(v)\,\mathrm{d}u\,\mathrm{d}v \leq Cb^{2\beta-\alpha}(1+r^{\alpha-2}), \quad (5.4)$$

*where the constant $C$ depends only on $\alpha$ and $\beta$.*

**Proof.** Observe that $F(b^{-1}u) = b^\beta F_1(u)$, where $F_1(u) = |u|^{-\beta}$ for $|u| < 1$ and $F_1(u) = |u|^{-1}$ for $|u| \geq 1$. Then, since

$$I_0(r) = \int_0^\infty \int_0^\infty |rb^{-1}u - b^{-1}v|^{\alpha-2} F(b^{-1}u)F(b^{-1}v)b^{-2}\,\mathrm{d}u\,\mathrm{d}v$$

$$= b^{2\beta-\alpha}\int_0^\infty \int_0^\infty |ru-v|^{\alpha-2} F_1(u)F_1(v)\,\mathrm{d}u\,\mathrm{d}v,$$

it is enough to show (5.4) with $b = 1$. Consider the integral $I_0(r)$ over four regions $(0,1) \times (0,1)$, $(1,\infty) \times (1,\infty)$, $(0,1) \times (1,\infty)$ and $(1,\infty) \times (0,1)$, and denote the corresponding integrals by $I_{0,1}(r)$, $I_{0,2}(r)$, $I_{0,3}(r)$ and $I_{0,4}(r)$, respectively. By making the changes of variables $v = ry, u = yw$ below, we obtain that

$$I_{0,1}(r) = \int_0^1 \int_0^1 |ru-v|^{\alpha-2} u^{-\beta} v^{-\beta}\,\mathrm{d}u\,\mathrm{d}v$$

$$= \int_0^{1/r} \mathrm{d}y \left(\int_0^{1/y} \mathrm{d}w |w-1|^{\alpha-2} w^{-\beta}\right) y^{\alpha-1-2\beta}\,r^{\alpha-1-\beta}.$$

First, suppose that $r \geq 1$. We examine $I_{0,1}(r)$ in three cases, depending on the behavior of the inner integral above.

*Case 1*: $\alpha - 2 - \beta + 1 = \alpha - 1 - \beta > 0$. In this case, since $r \geq 1$ and $\alpha < 2$, we have

$$I_{0,1}(r) \leq C \int_0^{1/r} y^{-(\alpha-1-\beta)} y^{\alpha-1-2\beta}\,\mathrm{d}y\,r^{\alpha-1-\beta} = C'r^{\alpha-2} \leq C'.$$

*Case 2*: $\alpha - 1 - \beta < 0$. In this case, since $r \geq 1$ and $\beta < 1$, we have

$$I_{0,1}(r) \leq C \int_0^{1/r} y^{\alpha-1-2\beta}\,\mathrm{d}y\,r^{\alpha-1-\beta} = C'r^{\beta-1} \leq C'.$$

*Case 3*: $\alpha - 1 - \beta = 0$. By using the integration by parts formula, we have

$$I_{0,1}(r) \leq C \int_0^{1/r} |\log y| y^{-\beta} \, dy \leq C'(r^{\beta-1} \log r + r^{\beta-1}) \leq C''.$$

This shows that $I_{0,1}(r) \leq C$ when $r \geq 1$. If $0 < r \leq 1$, then, by writing $I_{0,1}(r) = \int_0^1 \int_0^1 |u - r^{-1}v|^{\alpha-2} u^{-\beta} v^{-\beta} \, du \, dv \, r^{\alpha-2}$ and using symmetry, we can conclude, as above, that $I_{0,1}(r) \leq C r^{\alpha-2}$. Hence, $I_{0,1}(r) \leq C(1 + r^{\alpha-2})$ for any $r > 0$.

Turning to the integral $I_{0,2}(r)$, supposing first that $r \geq 1$ and making the change of variables $ru = y$ below, we obtain that

$$I_{0,2}(r) = \int_1^\infty \int_1^\infty |ru - v|^{\alpha-2} u^{-1} v^{-1} \, du \, dv = \int_r^\infty dy \int_1^\infty dv |y - v|^{\alpha-2} y^{-1} v^{-1} \leq C.$$

By using the symmetry argument above (see the case of $I_{0,1}(r)$), we can deduce that $I_{0,2}(r) \leq r^{\alpha-2}$ when $0 < r \leq 1$. Hence, $I_{0,2}(r) \leq C(1 + r^{\alpha-2})$ for any $r > 0$.

By making the changes of variables $v = ry, u = yw$, we obtain that

$$\begin{aligned} I_{0,3}(r) &= \int_0^1 du \int_1^\infty dv |ru - v|^{\alpha-2} u^{-\beta} v^{-1} \\ &= \int_{1/r}^\infty dy \left( \int_0^{1/y} dw |w - 1|^{\alpha-2} w^{-\beta} \right) y^{\alpha-2-\beta} \, r^{\alpha-2}. \end{aligned} \quad (5.5)$$

Consider first $r \geq 1$. In this case, we have

$$I_{0,3}(r) = \int_1^\infty dy \left( \int_0^{1/y} dw |w - 1|^{\alpha-2} w^{-\beta} \right) y^{\alpha-2-\beta} \, r^{\alpha-2}$$

$$+ \int_{1/r}^1 dy \left( \int_0^{1/y} dw |w - 1|^{\alpha-2} w^{-\beta} \right) y^{\alpha-2-\beta} \, r^{\alpha-2} =: I_{0,3}^1(r) + I_{0,3}^2(r).$$

For the integral $I_{0,3}^1(r)$, we have $I_{0,3}^1(r) = C \int_1^\infty y^{\beta-1} y^{\alpha-2-\beta} \, dy \, r^{\alpha-2} \leq C' r^{\alpha-2} \leq C'$ since $r \geq 1$. To bound the integral $I_{0,3}^2(r)$, we consider the following three cases.

*Case 1*: $\alpha - 1 - \beta > 0$. In this case, we have

$$I_{0,3}^2(r) \leq C \int_{1/r}^1 y^{-(\alpha-1-\beta)} y^{\alpha-2-\beta} \, dy \, r^{\alpha-2} = C'(\log r) r^{\alpha-2} \leq C'.$$

*Case 2*: $\alpha - 1 - \beta < 0$. In this case, we have

$$I_{0,3}^2(r) \leq C \int_{1/r}^1 y^{\alpha-2-\beta} \, dy \, r^{\alpha-2} = C' r^{\beta-1} \leq C'.$$

*Case 3*: $\alpha - 1 - \beta = 0$. In this case, we have

$$I_{0,3}^2(r) \leq C \int_{1/r}^1 |\log y| y^{-1} \, dy \, r^{\alpha-2} = C'(\log r)^2 r^{\alpha-2} \leq C''.$$

This shows that $I_{0,3}^2(r) \leq C$ when $r \geq 1$. Now, consider $0 < r \leq 1$. Since $r^{-1} \geq 1$, by using (5.5), we obtain that

$$I_{0,3}(r) \leq \int_1^\infty dy \left( \int_0^{1/y} dw |w-1|^{\alpha-2} w^{-\beta} \right) y^{\alpha-2-\beta} \, r^{\alpha-2} = Cr^{\alpha-2}.$$

Hence, $I_{0,3}(r) \leq C(1 + r^{\alpha-2})$ for any $r > 0$. It can be similarly shown that the same bound holds for $I_{0,4}(r)$. This proves the result (5.4). □

The following inequalities were used in the proofs of Lemmas 3.1, 3.3 and 5.1.

**Lemma 5.3.** *Let $x_1, x_2 \in \mathbb{R}$. When $\alpha \in (1,2)$, the following inequalities hold*

$$|x_1^{\langle \alpha-1 \rangle} - x_2^{\langle \alpha-1 \rangle}| \leq 2|x_2|^{\alpha-2}|x_1 - x_2| \qquad (x_2 \neq 0), \tag{5.6}$$

$$|x_1^{\langle \alpha-1 \rangle} - x_2^{\langle \alpha-1 \rangle}| \leq 2|x_1 - x_2|^{\alpha-1}. \tag{5.7}$$

*Moreover, when $\alpha \in (0,2)$,*

$$||x+y|^\alpha - |x|^\alpha - |y|^\alpha| \leq 2|xy|^{\alpha/2}. \tag{5.8}$$

**Proof.** The inequality (5.6) follows from the inequality $|(1+z)^{\langle q \rangle} - 1| \leq 2|z|$, with $q \in (0,1)$ and $z \in \mathbb{R}$, in Lemma 3.4 of Cioczek-Georges and Taqqu [7]. We provide here a geometric proof, which is best understood by graphing the function $p(x) = x^{\langle \alpha-1 \rangle}$, $x \in \mathbb{R}$. Suppose, without loss of generality, that $x_2 < 0$. Then, there exists an $x^*$ such that the slope of the line connecting the points $(x_2, x_2^{\langle \alpha-1 \rangle})$ and $(x^*, (x^*)^{\langle \alpha-1 \rangle})$ is the largest among the slopes connecting $(x_2, x_2^{\langle \alpha-1 \rangle})$ and $(x, x^{\langle \alpha-1 \rangle})$ with $x \in \mathbb{R}$, $x \neq x_2$. This $x^*$ satisfies $0 < x^* \leq |x_2|$. Then,

$$\frac{|x_1^{\langle \alpha-1 \rangle} - x_2^{\langle \alpha-1 \rangle}|}{|x_1 - x_2|} \leq \frac{|(x^*)^{\langle \alpha-1 \rangle} - x_2^{\langle \alpha-1 \rangle}|}{|x^* - x_2|} = \frac{|x^*|^{\alpha-1} + |x_2|^{\alpha-1}}{|x^*| + |x_2|} \leq \frac{2|x_2|^{\alpha-1}}{|x_2|} = 2|x_2|^{\alpha-2},$$

which is the inequality (5.6).

We now turn to the inequality (5.7). When $x_1, x_2 > 0$ or $x_1, x_2 < 0$, it follows from a sharper bound $|a_1^p - a_2^p| \leq |a_1 - a_2|^p$, valid for $a_1, a_2 > 0$ and $p \in (0,1)$. When $x_1 > 0, x_2 < 0$ or $x_1 < 0, x_2 > 0$, we have

$$|\langle x_1 \rangle^{\alpha-1} - \langle x_2 \rangle^{\alpha-1}| = |x_1|^{\alpha-1} + |x_2|^{\alpha-1} \leq 2(|x_1| + |x_2|)^{\alpha-1} = 2|x_1 - x_2|^{\alpha-1}.$$

Finally, the inequality (5.8) follows from $|x|^\alpha + |y|^\alpha - 2|xy|^{\alpha/2} = ||x|^{\alpha/2} - |y|^{\alpha/2}|^2 \leq |x+y|^{2\alpha/2} = |x+y|^\alpha$ and $|x+y|^\alpha = |x+y|^{2\alpha/2} = (|x|^2 + |y|^2 + 2|xy|)^{\alpha/2} \leq |x|^\alpha + |y|^\alpha + 2|xy|^{\alpha/2}$ since $0 < \alpha/2 < 1$. □

## 6. Application to central limit theorems

The following result now follows easily from the covariance bounds in Theorems 2.1–2.4 and the central limit theorems for bounded functionals of $S\alpha S$ moving averages considered in Pipiras and Taqqu [16] and Hsing [11].

**Theorem 6.1.** *Let $\alpha \in (1,2)$ and $\{\xi_n\}_{n\geq 0}$ be an $S\alpha S$ moving average defined by (2.3) with a function $a \in L^\alpha(\mathbb{R}, \mathrm{d}x)$, $a(x) = 0$ for $x < 0$. Suppose that $K$ is a function satisfying Assumptions* (K1)–(K3). *If*

$$\sum_{m=1}^{\infty} \left(\int_{m-1}^{m} |a(x)|^\alpha \,\mathrm{d}x\right)^{1/2} < \infty, \qquad \sum_{n=1}^{\infty} [\xi_0, \xi_n]_1 < \infty, \tag{6.1}$$

*then*

$$N^{-1/2} S_N = N^{-1/2} \sum_{n=1}^{N} (K(\xi_n) - EK(\xi_n)) \xrightarrow{d} \mathcal{N}(0, \sigma^2) \tag{6.2}$$

*as $N \to \infty$, where*

$$\sigma^2 = \lim_{N\to\infty} N^{-1} \mathrm{Var}(S_N) = \mathrm{Var}(K(\xi_0)) + 2\sum_{n=1}^{\infty} \mathrm{Cov}(K(\xi_0), K(\xi_n)). \tag{6.3}$$

*The limit in* (6.3) *exists and is finite and the series in* (6.3) *converges absolutely.*

**Remark 6.1.** The assumption of Theorem 6.1 that moving averages are causal (i.e., $a(x) = 0$ for $x < 0$) is also present in the available central limit theorems for bounded functions $K$ (see Hsing [11], Pipiras and Taqqu [16]).

The easiest way to verify condition (6.1) is to use, when possible, the stronger condition stated in the following lemma. Its proof is elementary and is given below. We use the notation $[\xi_0, \xi_n]_{\gamma,\delta} = \int_\mathbb{R} |a(-x)|^\gamma |a(n-x)|^\delta \,\mathrm{d}x$.

**Lemma 6.1.** *Let $\gamma > 0$, $\delta > 0$ be such that $\gamma + \delta = \alpha$ with $\alpha > 0$. Also, let $a \in L^\alpha(\mathbb{R}, \mathrm{d}x)$. Then,*

$$\sum_{m=-\infty}^{\infty} \left(\int_{m-1}^{m} |a(x)|^\alpha \,\mathrm{d}x\right)^{\min\{\gamma,\delta\}/\alpha} < \infty \quad \Longrightarrow \quad \sum_{n=-\infty}^{\infty} [\xi_0, \xi_n]_{\gamma,\delta} < \infty. \tag{6.4}$$

*In particular, when $\alpha \in (1,2)$, we have*

$$\sum_{m=-\infty}^{\infty} \left(\int_{m-1}^{m} |a(x)|^\alpha \,\mathrm{d}x\right)^{(\alpha-1)/\alpha} < \infty \quad \Longrightarrow \quad \sum_{n=1}^{\infty} [\xi_0, \xi_n]_1 < \infty, \tag{6.5}$$

$$\sum_{m=-\infty}^{\infty} \left( \int_{m-1}^{m} |a(x)|^\alpha \, dx \right)^{1/2} < \infty \implies \sum_{n=1}^{\infty} [\xi_0, \xi_n]_2 < \infty. \qquad (6.6)$$

**Proof.** By setting $I = \sum_{n=1}^{\infty} [\xi_0, \xi_n]_{\gamma,\delta}$, the result (6.4) follows from Hölder's inequality,

$$I = \sum_{m,n \in \mathbb{Z}} \int_{m-1}^{m} |a(z)|^\gamma |a(n+z)|^\delta \, dz$$

$$\leq \sum_{m,n \in \mathbb{Z}} \left( \int_{m-1}^{m} |a(z)|^\alpha \, dz \right)^{\gamma/\alpha} \left( \int_{m-1}^{m} |a(n+z)|^\alpha \, dz \right)^{\delta/\alpha}$$

$$= \sum_{m \in \mathbb{Z}} \left( \int_{m-1}^{m} |a(z)|^\alpha \, dz \right)^{\gamma/\alpha} \sum_{n \in \mathbb{Z}} \left( \int_{n+m-1}^{n+m} |a(z)|^\alpha \, dz \right)^{\delta/\alpha},$$

making the substitution $n \to n - m$ in the last summation. $\square$

**Example 6.1.** Let $\alpha \in (1,2)$ and $a: \mathbb{R} \mapsto \mathbb{R}$ be a bounded function such that $a(x) = 0$ for $x < 0$ and $|a(x)| \leq Cx^{-p}$ for $x > 0$ and some $p > 1/\alpha$. Then, for $m \geq 1$ and some constant $C$,

$$\int_{m-1}^{m} |a(x)|^\alpha \, dx \leq C m^{-p\alpha}.$$

Since $0 < (\alpha - 1)/\alpha < 1/2$, the conditions on the left-hand sides of (6.6) and (6.5) hold as long as $\sum_{m=1}^{\infty} |m|^{-p(\alpha-1)} < \infty$ or, equivalently, $p > 1/(\alpha - 1)$. Then, by (6.5), the function $a$ satisfies condition (6.1) of Theorem 6.1 for $p > 1/(\alpha - 1)$.

**Proof of Theorem 6.1.** Suppose, for simplicity, that $p = 1$ in Assumption (K3) and set $K_b(x) = K(x) 1_{\{|x| > b\}}$ and $K_{1,b}(x) = K(x) 1_{\{|x - x_1| < 1/b\}}$. Then, the function $L_b(x) = K(x) - K_b(x) - K_{1,b}(x)$ is bounded. Denote the partial sum in (6.2) by $S_{b,N}$ when $K$ is replaced by $L_b$. To show (6.2), by Theorem 4.2 of Billingsley [6], it is enough to prove that:

(i) $N^{-1/2} S_{b,N} \xrightarrow{d} \mathcal{N}(0, \sigma_b^2)$ as $N \to \infty$;
(ii) $\sigma_b^2 \to \sigma^2$ as $b \to \infty$;
(iii) $\limsup_{b \to \infty} \limsup_{N \to \infty} N^{-1} \text{Var}(S_N - S_{b,N}) = 0$.

Part (i) follows from Theorem 2.1 of Pipiras and Taqqu [16] since the first condition in (6.1) holds and the function $L_b$ is bounded. Moreover,

$$\sigma_b^2 = \lim_{N \to \infty} N^{-1} \text{Var}(S_{b,N}). \qquad (6.7)$$

Part (ii) will follow from (iii). Indeed, since

$$|\text{Var}^{1/2}(S_{b_1,N}) - \text{Var}^{1/2}(S_{b_2,N})| \leq \text{Var}^{1/2}(S_{b_1,N} - S_N) + \text{Var}^{1/2}(S_{b_2,N} - S_N),$$

*Bounds for the covariance of functions of infinite variance stable random variables* 1115

if (iii) holds, then $\sigma_b$ converges to some $\sigma$ by the Cauchy criterion. To prove (iii), observe that

$$S_N - S_{b,N} = \sum_{n=1}^{N}(K_b(\xi_n) - EK_b(\xi_n)) + \sum_{n=1}^{N}(K_{1,b}(\xi_n) - EK_{1,b}(\xi_n))$$

and hence

$$N^{-1}\operatorname{Var}(S_N - S_{b,N}) \leq 2N^{-1}\operatorname{Var}\left(\sum_{n=1}^{N}(K_b(\xi_n) - EK_b(\xi_n))\right)$$

$$+ 2N^{-1}\operatorname{Var}\left(\sum_{n=1}^{N}(K_{1,b}(\xi_n) - EK_{1,b}(\xi_n))\right)$$

$$= 2\operatorname{Var}(K_b(\xi_0)) + 4\sum_{n=1}^{N-1}\left(1 - \frac{n}{N}\right)\operatorname{Cov}(K_b(\xi_0), K_b(\xi_n))$$

$$+ 2\operatorname{Var}(K_{1,b}(\xi_0)) + 4\sum_{n=1}^{N-1}\left(1 - \frac{n}{N}\right)\operatorname{Cov}(K_{1,b}(\xi_0), K_{1,b}(\xi_n)),$$

which implies that

$$\limsup_{N\to\infty} N^{-1}\operatorname{Var}(S_N - S_{b,N}) \leq 2\operatorname{Var}(K_b(\xi_0)) + 4\sum_{n=1}^{\infty}|\operatorname{Cov}(K_b(\xi_0), K_b(\xi_n))|$$

$$+ 2\operatorname{Var}(K_{1,b}(\xi_0)) + 4\sum_{n=1}^{\infty}|\operatorname{Cov}(K_{1,b}(\xi_0), K_{1,b}(\xi_n))|.$$

One can bound these sums by using Theorems 2.2 and 2.1 and Proposition 2.2. Sums involving the bounds $[\xi_0, \xi_n]_1$ converge by the second relation in (6.1) and the sum involving the bound $[\xi_0, \xi_n]_2$ converges by the first relation in (6.1) and Lemma 6.1 above. Part (iii) then follows by letting $b \to \infty$.

To show that $\sigma^2$ is the asymptotic variance of $N^{-1/2}S_N$ in (6.3), first observe that

$$|N^{-1/2}\operatorname{Var}^{1/2}(S_N) - \sigma| \leq |N^{-1/2}\operatorname{Var}^{1/2}(S_N - S_{b,N})| + |N^{-1/2}\operatorname{Var}^{1/2}(S_{b,N}) - \sigma|.$$

Then, the first relation in (6.3) follows by taking $\limsup_{b\to\infty} \limsup_{N\to\infty}$ in the above inequality and using (6.7) and parts (ii) and (iii) above. Finally, we need to show that the series in (6.3) converges absolutely. By Theorem 2.4 and Proposition 2.2, we get, for sufficiently large $n$,

$$|\operatorname{Cov}(K(\xi_0), K(\xi_n))| \leq C([\xi_0, \xi_n]_1 + [\xi_0, \xi_n]_2).$$

The absolute convergence then follows by using, as before, the relations (6.1) and Lemma 6.1. □

**Corollary 6.1.** *Theorem 6.1 is also valid for $S\alpha S$ moving averages $\{\xi_n\}_{n\geq 0}$ defined in (2.3) by using function $a$ such that $a(x) = 0$ when $x < x_0$, for some $x_0 < 0$.*

**Proof.** Observe that the kernel $a(n-x)$ in (2.3) can be replaced by $a(n-(x-x_0)) =: \widetilde{a}(n-x)$ without changing the distribution of $\{\xi_n\}_{n\geq 0}$. The function $\widetilde{a}$ is such that $\widetilde{a}(x) = 0$ for $x < 0$ and, if the function $a$ satisfies (6.1), then so does the function $\widetilde{a}$. Then, by applying Theorem 6.1 to the moving average defined via the function $\widetilde{a}$, we conclude that (6.2) and (6.3) hold. □

*Remark 6.2.* A related paper of Wu [23] also contains a central limit theorem for infinite variance causal moving averages (not necessarily stable) of the "discrete" form $X_n = \sum_{k=0}^{\infty} a_k \varepsilon_{n-k}$ and possibly unbounded functions $K$; see Theorem 1(a) in Wu [23]. The method of proof in Wu [23] is based on general central limit theorems for Markov chains developed by Michael Woodroofe and is quite different from the one presented here. The following observations shed light on the relationship between the result of Wu [23] and Theorem 6.1 above.

First, there are unbounded functions $K$ that satisfy the assumptions of Theorem 6.1 but not those of Theorem 1(a) of Wu [23]. An example is the function $K(x) = \log_2 |x| - E \log_2 |X_0|$, of particular interest in the wavelet-based application discussed below. Theorem 1(a) in Wu [23] involves the function $L_{K_n}(x)$, where $K_n(w) = EK(w + \overline{X}_{n,1})$, $\overline{X}_{n,1} = a_0 \varepsilon_n + a_1 \varepsilon_{n-1} + \cdots + a_{n-1} \varepsilon_1$ and

$$L_f(x) = \sup_{y:|y-x|<1} \frac{|f(y) - f(x)|}{|y-x|}$$

is a local Lipschitz constant. For the logarithm function $K$ above,

$$L_{K_n}(x) \geq \frac{1}{\ln 2} E \frac{1}{|x + \overline{X}_{n,1}|} = \infty$$

and hence Theorem 1(a) of Wu [23] cannot be applied.

Second, consider, for example, the power function $K(x) = |x|^\beta - E|X_0|^\beta$ with $\beta < \alpha/2$, also of interest in wavelet-based application. Theorem 1(a) of Wu [23] applies to this power function under suitable conditions. We next examine how these conditions compare to those in Theorem 6.1 above. First, we can verify that the condition on $L_{K_n}$ in Wu [23] is satisfied when $(\beta - 1)2q/(q-1) > -1$ or $1/q < 2\beta - 1$. Observe that since Wu requires $q > 1$, we obtain $2\beta - 1 > 0$ and hence $\alpha > 1$, as in Theorem 6.1 above. Second, the condition on $K_n$ in Wu [23] can be seen to be verified when $\beta 2q/(q-1) < \alpha$ or $1/q < (\alpha - 2\beta)/\alpha$. Comparing the conditions on $L_{K_n}$ and $K_n$, observe that they are verified if $1/q < 2\beta - 1$ when $1/2 < \beta < \alpha/(\alpha+1)$, and if $1/q < (\alpha - 2\beta)/\alpha$ when $\alpha/(\alpha+1) < \beta < \alpha/2$.

Now, suppose that the moving average coefficients $a_n$ satisfy $|a_n| \sim Cn^{-p}$ as $n \to \infty$, where $p > 1/\alpha$. The condition on $a_n$ in Wu [23] requires that $\sum_n |a_n|^{\alpha/2q} < \infty$ or $1 < p\alpha/2q$. From the bounds on $1/q$ above, it is therefore necessary to have $1 < p\alpha(2\beta - 1)/2$ when $1/2 < \beta < \alpha/(\alpha+1)$, and $1 < p(\alpha - 2\beta)/2$ when $\alpha/(\alpha+1) < \beta < \alpha/2$. On the other hand, by Example 6.1 above, the conditions on $a_n$ of Theorem 6.1 are satisfied when $1 < p(\alpha - 1)$. When $\alpha/(\alpha + 1) < \beta < \alpha/2$, observe that $\alpha - 2\beta < 2(\alpha - 1)$. Hence, Theorem 6.1 is stronger than Theorem 1(a) in Wu [23]. When $1/2 < \beta < \alpha/(\alpha + 1)$, observe that $\alpha(2\beta - 1) < 2(\alpha - 1)$ and hence Theorem 6.1 remains stronger. Combining these observations, Theorem 6.1 is stronger than Theorem 1(a) of Wu [23] for the above power function. Presently, we are not aware of any unbounded functions where Theorem 1(a) of Wu [23] performs better.

Theorem 6.1 has the following multivariate extension which is used in the next section. The proof of the extension is analogous to that of Theorem 6.1 and is omitted, except for a supplementary result on the form of a limit covariance. Let $c_j > 0$ be positive real numbers for $j = 1, \ldots, J$. Consider the $S\alpha S$ moving average sequences $\{\xi_{j,n}\}_{n \geq 0}$, $j = 1, \ldots, J$, given by

$$\xi_{j,n} = \int_{\mathbb{R}} a_j(n - c_j x) M(\mathrm{d}x), \tag{6.8}$$

where $a_j \in L^\alpha(\mathbb{R}, \mathrm{d}x)$ and $a_j(x) = 0$ for $x < 0$ and $M$ is an $S\alpha S$ random measure with the Lebesgue control measure on $\mathbb{R}$. Fix $n_j$, $j = 1, \ldots, J$, and let $N_j$ be positive integers such that, as $N \to \infty$,

$$N_j \sim \frac{N}{n_j}. \tag{6.9}$$

**Theorem 6.2.** *Let $\alpha \in (1, 2)$ and $\{\xi_{j,n}\}_{n \geq 0}$, $j = 1, \ldots, J$, be $S\alpha S$ moving averages defined by (6.8) with $a_j$ such that $a_j(x) = 0$ when $x < x_0$, for some fixed $x_0$. Suppose that, for each $j = 1, \ldots, J$, the function $K_j$ and the kernel $a_j$ in (6.8) satisfy conditions (K1)–(K3) and (6.1), respectively. Then, as $N \to \infty$,*

$$(N_j^{-1/2} S_{j,N_j})_{j=1}^J = \left( N_j^{-1/2} \sum_{n=1}^{N_j} (K_j(\xi_{j,n}) - EK_j(\xi_{j,n})) \right)_{j=1}^J \xrightarrow{d} \mathcal{N}(\mathbf{0}, \boldsymbol{\sigma}), \tag{6.10}$$

*where $\boldsymbol{\sigma} = (\sigma_{jk})_{j,k=1,\ldots,J}$ with*

$$\sigma_{jk} = \lim_{N \to \infty} E \frac{S_{j,N_j}}{N_j^{1/2}} \frac{S_{k,N_k}}{N_k^{1/2}}, \tag{6.11}$$

*which exists and is finite.*

*If, in addition, $c_j = 2^{-j}$, $n_j = 2^j$ and, for all $k > j$, $\sum_{n=1}^\infty [\xi_{j,n}, \xi_{k,0}]_1 < \infty$, then the asymptotic covariance (6.11) is given by*

$$\sigma_{jk} = 2^{(j-k)/2} \sum_{n=-\infty}^\infty \mathrm{Cov}(K_j(\xi_{j,n}), K_k(\xi_{k,0})), \qquad k \geq j, \tag{6.12}$$

*where the series in (6.12) converges absolutely.*

**Proof.** We prove the series representation (6.12). Recall that we now suppose $c_j = 2^{-j}$ and $n_j = 2^j$. We shall also assume, for simplicity, that $n_j N_j = N$ in (6.9). (The general case can be proven in a similar way.) Then, $N_j = 2^{k-j} N_k$ and we have, for $k \geq j$,

$$
\begin{aligned}
ES_{j,N_j} S_{k,N_k} &= \sum_{n_1=1}^{N_j} \sum_{n_2=1}^{N_k} \mathrm{Cov}(K_j(\xi_{j,n_1}), K_k(\xi_{k,n_2})) \\
&= \sum_{n_1=1}^{2^{k-j} N_k} \sum_{n_2=1}^{N_k} \mathrm{Cov}(K_j(\xi_{j,n_1-2^{k-j} n_2}), K_k(\xi_{k,0})) \\
&= \sum_{p=0}^{2^{k-j}-1} \sum_{n_1=1}^{N_k} \sum_{n_2=1}^{N_k} \mathrm{Cov}(K_j(\xi_{j,2^{k-j} n_1 - p - 2^{k-j} n_2}), K_k(\xi_{k,0})) \\
&= \sum_{p=0}^{2^{k-j}-1} \sum_{n=-N_k}^{N_k} (N_k - |n|) \mathrm{Cov}(K_j(\xi_{j,2^{k-j} n - p}), K_k(\xi_{k,0})).
\end{aligned}
$$

After dividing by $N_j^{1/2} N_k^{1/2} = 2^{(k-j)/2} N_k$ and letting $N \to \infty$, we then obtain by (6.11), that

$$
\begin{aligned}
\sigma_{jk} &= 2^{(j-k)/2} \sum_{p=0}^{2^{k-j}-1} \sum_{n=-\infty}^{\infty} \mathrm{Cov}(K_j(\xi_{j,2^{k-j} n - p}), K_k(\xi_{k,0})) \\
&= 2^{(j-k)/2} \sum_{n=-\infty}^{\infty} \mathrm{Cov}(K_j(\xi_{j,n}), K_k(\xi_{k,0})),
\end{aligned}
$$

provided the series is absolutely convergent. This can be established by first using a result for $\xi_{j,n}$, analogous to Proposition 2.2 and Theorem 2.4, to conclude that, for sufficiently large $n$,

$$|\mathrm{Cov}(K_j(\xi_{j,n}), K_k(\xi_{k,0}))| \leq C([\xi_{j,n}, \xi_{k,0}]_1 + [\xi_{j,n}, \xi_{k,0}]_2). \tag{6.13}$$

When $k = j$, the convergence of the series with the terms on the right-hand side of (6.13) follows from the fact that the $a_j$'s satisfy (6.1). When $k > j$, we use the assumption $\sum_n [\xi_{j,n}, \xi_{k,0}]_1 < \infty$, together with the observation that the first relation in (6.1) with $a = a_j$ and $a = a_k$ implies $\sum_n [\xi_{j,n}, \xi_{k,0}]_2 < \infty$, which can be proven as was (6.6) in Lemma 6.1. □

## 7. Asymptotic normality of wavelet-based estimators

We here apply Theorem 6.2 to prove the asymptotic normality of wavelet-based estimators of the self-similarity parameter in linear fractional stable motion (LFSM). LFSM is an $S\alpha S$ self-similar stationary-increments process with the integral representation

$$X(t) = \int_{\mathbb{R}} \{(t-u)_+^\kappa - (-u)_+^\kappa\} M(\mathrm{d}u), \qquad t \in \mathbb{R}, \tag{7.1}$$

where

$$\kappa = H - 1/\alpha, \tag{7.2}$$

$H \in (0,1)$ is the self-similarity parameter, $\alpha \in (0,2)$ is the stability parameter and $M$ is an $S\alpha S$ random measure with the Lebesgue control measure on $\mathbb{R}$. Self-similarity means that for all $c > 0$, the processes $X(ct)$ and $c^H X(t)$ have the same finite-dimensional distributions. LFSM is an infinite variance counterpart of fractional Brownian motion, which is the only Gaussian self-similar process with stationary increments. It is often taken as a representative process for self-similar processes with stationary increments having infinite variance. (The more general double-sided LFSM cannot be considered for Theorem 7.1 below because of the causality assumption discussed in Remark 6.1.) For more information on LFSM, see Section 7 of Samorodnitsky and Taqqu [19].

Wavelets have already proven useful for estimating the self-similarity parameter of fractional Brownian motion and a related long-memory parameter in finite variance long-memory time series (Veitch and Abry [22], Bardet [5] and Abry, Flandrin, Taqqu and Veitch [2]). This work was motivated to a great extent by applications to data traffic in communication networks. Following this line of research, several authors have suggested using wavelets to estimate the self-similarity parameter $H$ in LFSM as well. From a mathematical point of view, it is interesting to see to what extent existing wavelet methodology applies to this extended class of processes. There is indeed evidence that some teletraffic data deviates from Gaussianity and has heavy tail, characteristic of processes with infinite variance.

The discrete wavelet transform of LFSM (or of other deterministic or stochastic functions) is a sequence of *discrete wavelet transform* coefficients $\{d_{j,k}\}_{j,k \in \mathbb{Z}}$, defined by

$$d_{j,k} = \int_{\mathbb{R}} X(t) \psi_{j,k}(t) \, \mathrm{d}t = \int_{\mathbb{R}} X(t) 2^{-j/2} \psi(2^{-j} t - k) \, \mathrm{d}t. \tag{7.3}$$

Here, $\psi : \mathbb{R} \mapsto \mathbb{R}$ is the so-called *wavelet function* which has $Q \geq 1$ zero moments, that is,

$$\int_{\mathbb{R}} \psi(t) t^m \, \mathrm{d}t = 0, \qquad m = 0, \ldots, Q-1 \quad \text{and} \quad \int_{\mathbb{R}} \psi(t) t^Q \, \mathrm{d}t \neq 0. \tag{7.4}$$

It was shown by Delbeke and Abry [8, 9] and Pesquet-Popescu [15] that the discrete wavelet transform coefficients of LFSM are well defined for a bounded, compactly supported wavelet $\psi$ when $H - 1/\alpha > -1$ and that the following Fubini-type result holds for $j, k \in \mathbb{Z}$:

$$d_{j,k} = \int_{\mathbb{R}} \left( \int_{\mathbb{R}} (t-u)_+^\kappa \psi_{j,k}(t) \, \mathrm{d}t \right) M(\mathrm{d}u) = \int_{\mathbb{R}} 2^{j(\kappa+1/2)} h(k - 2^{-j}u) M(\mathrm{d}u), \tag{7.5}$$

where

$$h(u) = \int_{\mathbb{R}} (s+u)_+^\kappa \psi(s) \, ds, \qquad u \in \mathbb{R}. \tag{7.6}$$

Moreover, it is easy to see that for each $j \in \mathbb{Z}$, the sequence $\{d_{j,k}\}_{k \in \mathbb{Z}}$ is stationary, that is, for all $l \in \mathbb{Z}$, $\{d_{j,k+l}\} =_d \{d_{j,k}\}$, and that the following scaling relation holds:

$$\{d_{j,k}\}_{k \in \mathbb{Z}} \stackrel{d}{=} \{2^{j(H+1/2)} d_{0,k}\}_{k \in \mathbb{Z}}. \tag{7.7}$$

Let us also note that, in practice, the wavelet coefficients are computed (rather, approximated) by using fast pyramidal Mallat-type algorithms (see, e.g., Mallat [14]).

A wavelet-based estimator of $H$ can be defined by

$$\widehat{H} = \sum_j w_j \frac{1}{N_j} \sum_{n=1}^{N_j} \log_2 |d_{j,n}| - \frac{1}{2}. \tag{7.8}$$

The summation here is over some finite number of consecutive indices $j$, called *octaves* or *scales*, and each integer $N_j$ corresponds to the number of available wavelet coefficients of LFSM over the time interval $[0, N]$ at scale $j$. Roughly speaking, we have $N_j = 2^{-j} N$ (up to border effects) and thus $N_j = N_j(N) \to \infty$ as $N \to \infty$. Finally, the $w_j$'s satisfy the relations

$$\sum_j w_j = 0, \qquad \sum_j j w_j = 1 \tag{7.9}$$

and can be viewed as weights for a linear least square estimation of a slope in $(j, Y_j)$, where $Y_j$ is the term multiplying $w_j$ in (7.8). While the wavelet coefficients $d_{j,k}$ have infinite variance, $\log_2 |d_{j,k}|$ has finite variance, as do the estimators $\widehat{H}$. Theorem 6.2 implies that this estimator is asymptotically normal.

**Theorem 7.1.** *Let $\widehat{H}$ be the self-similarity parameter estimator in LFSM defined by (7.8) for bounded and compactly supported wavelets. Suppose that $\alpha \in (1, 2)$ and that*

$$Q - H > \frac{1}{\alpha(\alpha-1)}. \tag{7.10}$$

*Then,*

$$N^{1/2}(\widehat{H} - H) \stackrel{d}{\to} \mathcal{N}(0, \sigma^2) \tag{7.11}$$

*as $N \to \infty$, where*

$$\sigma^2 = \sum_{j,k} w_j w_k 2^{j/2} 2^{k/2} \sigma_{j,k},$$

*with*

$$\sigma_{j,k} = 2^{(j-k)/2} \sum_{n=-\infty}^{\infty} \mathrm{Cov}(\log_2 |d_{j,n}|, \log_2 |d_{k,0}|) \tag{7.12}$$

if $k \geq j$, and $\sigma_{j,k} = \sigma_{k,j}$ if $j > k$. We also have

$$\sigma^2 = \lim_{N \to \infty} NE(\widehat{H} - H)^2. \tag{7.13}$$

**Remark 7.1.** Because convergence in (7.11) is to a normal law, a rate of $N^{1/2}$ is expected. In Kokoszka and Taqqu [13], the Whittle method is used to estimate the long-range dependence parameter in stable FARIMA time series. A faster rate of essentially $N^{1/\alpha}$ is obtained, but also a limit which has infinite variance and is hence more spread out.

**Proof of Theorem 7.1.** Without loss of generality, we suppose that the sum in (7.8) is over octaves $j = 1, \ldots, J$ for some fixed $J$. By using (7.7), we can write

$$N^{1/2}(\widehat{H} - H) = \sum_j w_j 2^{j/2} \frac{1}{N_j^{1/2}} \sum_{n=1}^{N_j} (\log_2 |d_{j,n}| - E \log_2 |d_{j,n}|) = \mathbf{w}\mathbf{Y}^t, \tag{7.14}$$

where $t$ stands for "transposed," $\mathbf{w} = (2^{1/2}w_1, \ldots, 2^{J/2}w_J)$ and $\mathbf{Y} = (Y(1, N_1), \ldots, Y(J, N_J))$, with

$$Y(j, N_j) = N_j^{-1/2} \sum_{n=1}^{N_j} (\log_2 |d_{j,n}| - E \log_2 |d_{j,n}|).$$

Therefore, to prove that $\widehat{H}$ is asymptotically normal, it is enough to establish the asymptotic normality of the vector $\mathbf{Y}$. Observe that, by (7.5), this vector can be expressed as that in (6.10) and (6.8) with $c_j = 2^{-j}$ and $a_j(x) = 2^{j(\kappa+1/2)}h(x)$, where $j = 1, \ldots, J$. Moreover, since the wavelet $\psi$ has compact support, there exists an $x_0 \in \mathbb{R}$ such that $h(x) = 0$ for $x < x_0$ and thus $a_j(x) = 0$ for $x < x_0$ as well. Then, by Theorem 6.2, we get asymptotic normality, provided the function $\log_2 |x|$ satisfies (K1)–(K3) and the function $a(x) = h(x)$ satisfies (6.1). It is easy to verify that $K(x) = \log_2 |x|$ satisfies (K1)–(K3). The function $h$ satisfies (6.1) by using Lemma 7.1 below and Example 6.1, provided that $Q - \kappa > 1/(\alpha - 1)$ or, equivalently, $Q - H > 1/(\alpha(\alpha - 1))$.

Let us now show that $\sigma^2$ can be expressed as in the theorem. By using (7.14) and Theorem 6.2, it is enough to show that under the assumption (7.10), $\sum_{n=1}^{\infty} [d_{j,n}, d_{k,0}]_1 < \infty$ for all $k > j$. This can be established as in Example 6.1 by using $a_j(x) = 2^{j(\kappa+1/2)}h(x)$, Lemma 7.1 below and the implication (6.5). $\square$

The following result was used in the proof of Theorem 7.1 above.

**Lemma 7.1.** *Let $h$ be the function defined by (7.6), where the wavelet $\psi$ has compact support, is bounded and has $Q \geq 1$ zero moments. Then, for sufficiently large $u$, $|h(u)| \leq Cu^{\kappa-Q}$, where $C$ is some constant.*

**Proof.** By applying Taylor's formula to the function $f_u(s) = (s+u)_+^\kappa$, which is infinitely differentiable on any interval $s \in [-M, M]$, $M > 0$, for sufficiently large $u$, we have that

$$(s+u)_+^\kappa = u_+^\kappa + \cdots + \frac{\kappa(\kappa-1)\cdots(\kappa-Q+2)}{(Q-1)!} u_+^{\kappa-Q+1} s^{Q-1}$$
$$+ \frac{\kappa(\kappa-1)\cdots(\kappa-Q+1)}{Q!}(s_0+u)_+^{\kappa-Q} s^Q,$$

where $s_0 \in [-M, M]$. If the wavelet $\psi$ has compact support, we then deduce the result by using (7.4). $\square$

Note that the greater the value of $Q$ in Lemma 7.1, the faster $|h(u)|$ decreases as $u$ grows. Thus, in view of (7.5), when $Q$ increases, the sequence $\{d_{j,k}\}_{k\in\mathbb{Z}}$ becomes almost independent. This is one of the main advantages of working with the wavelet coefficients rather than with LFSM itself. A more comprehensive and applied study of wavelet-based estimators of the self-similarity parameter in LFSM (and further results) can be found in Abry, Delbeke and Flandrin [1], Abry, Pesquet-Popescu and Taqqu [4], Abry, Fladrin, Taqqu and Veitch [3], Stoev, Pipiras and Taqqu [21] and Stoev and Taqqu [20].

*Remark 7.2.* An alternative estimator for the self-similarity parameter in LFSM can be defined as

$$\widehat{H}^* = \frac{1}{\beta} \sum_j w_j \log_2\left(\frac{1}{N_j}\sum_{n=1}^{N_j} |d_{j,n}|^\beta\right) - \frac{1}{2}, \qquad -1 < \beta < \alpha/2. \tag{7.15}$$

The asymptotic normality result for $\widehat{H}^*$ can be established as in Theorem 7.1 with

$$\sigma_{j,k} = 2^{(j-k)/2} \sum_{n=-\infty}^{\infty} \frac{\mathrm{Cov}(|d_{j,n}|^\beta, |d_{k,0}|^\beta)}{(\ln 2)^2 E|d_{j,0}|^\beta E|d_{k,0}|^\beta} \tag{7.16}$$

if $k \geq j$, and $\sigma_{j,k} = \sigma_{k,j}$ if $j > k$. However, because of the presence of the function $\log_2$ in (7.15), we cannot deduce that the relation (7.13) holds, namely that $\sigma^2$ is the asymptotic (normalized) variance of $\widehat{H}^*$.

*Remark 7.3.* Another interesting question not addressed here is that of joint estimation of $H$ and $\sigma$. We plan to pursue this question in future work.

## Acknowledgements

The authors would like to thank the associate editor and the two referees for many useful comments and suggestions. Thanks also to Michael Woodroofe for bringing the paper of Wu [23] to our attention.

The work was partially supported by the NSF Grants DMS-05-05628 at the University of North Carolina and DMS-0505747 at Boston University. The CNRS support for the visiting researcher position at ENS Lyon for V. Pipiras is also gratefully acknowledged.

[18] Rosiński, J. (2006). Minimal integral representations of stable processes. *Probab. Math. Statist.* **26** 121–142. MR2301892

[19] Samorodnitsky, G. and Taqqu, M.S. (1994). *Stable Non-Gaussian Processes: Stochastic Models with Infinite Variance.* London: Chapman and Hall. MR1280932

[20] Stoev, S. and Taqqu, M.S. (2005). Asymptotic self-similarity and wavelet estimation for long-range dependent FARIMA time series with stable innovations. *J. Time Ser. Anal.* **26** 211–249. MR2122896

[21] Stoev, S., Pipiras, V. and Taqqu, M.S. (2002). Estimation of the self-similarity parameter in linear fractional stable motion. *Signal Processing* **82** 873–1901.

[22] Veitch, D. and Abry, P. (1999). A wavelet-based joint estimator of the parameters of long-range dependence. *IEEE Trans. Inform. Theory* **45** 878–897. MR1682517

[23] Wu, W.B. (2003). Additive functionals of infinite-variance moving averages. *Statist. Sinica* **13** 1259–1267. MR2026072